\documentclass[preprint,10pt]{elsarticle}

\makeatletter
\def\ps@pprintTitle{%
 \let\@oddhead\@empty
 \let\@evenhead\@empty
 \def\@oddfoot{}%
 \let\@evenfoot\@oddfoot}
\makeatother

\usepackage{amssymb}
\usepackage[utf8]{inputenc}
\usepackage{amsmath}
\usepackage{amsfonts}
\usepackage{enumitem}  
\usepackage[ruled,vlined,linesnumbered,resetcount]{algorithm2e}
\usepackage{tikz}
\usepackage{caption}
\usepackage{subcaption}
\usepackage{float}

\usepackage{graphicx}
\usepackage[normalem]{ulem}
\usepackage{xcolor}
\usepackage{bm}

\newtheorem{remark}{Remark}

\newcommand{\myfrac}[2]{\displaystyle{\frac{#1}{#2}}}
\newcommand{\somme}[3]{\displaystyle{\overset{#3}{\underset{#1=#2}{\sum}}}}
\newcommand{\deriveedot}[1]{\displaystyle{\dot{#1}}}
\newcommand{\derivee}[2]{\displaystyle{\frac{d #1}{d #2}}}
\newcommand{\deriveepartielle}[2]{\displaystyle{\partial_{#2} #1}}

\newcommand{\ctrlvar}{\gamma}

\newcommand{\statevecU}{\vect{u}}

\newcommand{\stateU}{u}

\newcommand{\vect}[1]{ {\bm{#1}}}

\newcommand{\InterpSbase}{{\Phi_{_{\untrainedparam}}}}
\newcommand{\InterpTbase}{{\Lambda_{_{\untrainedparam}}}}

\newcommand{\Smode}[2]{\displaystyle{\varphi_{#1}^{#2}}}

\newcommand{\Tmode}[2]{{\alpha_{#1}^{#2}}}

\newcommand{\Sbase}[1]{{\Phi_{_{\trainingparam{#1}}}}}

\newcommand{\Xbasis}[1]{{\Phi_{_{#1}}}}
\newcommand{\XbasisEnd}[1]{{\Psi_{_{#1}}}}
\newcommand{\Tbasis}[1]{{\Lambda_{_{#1}}}}

\newcommand{\Tbase}[1]{{\Lambda_{_{\trainingparam{#1}}}}}
\newcommand{\trainedSingVal}[1]{{\Sigma_{_{\trainingparam{#1}}}}}
\newcommand{\InterpSingVal}{{\Sigma_{_{\untrainedparam}}}}

\newcommand{\integOmega}[1]{\displaystyle{\int_{\Omega}#1 \,dx}}

\newcommand{\nbrModes}{q}
\newcommand{\nbrModesU}{q_{u}}
\newcommand{\nbrModesP}{q_{p}}

\newcommand{\nbrParam}{N_{p}}
\newcommand{\nbrSnap}{N_{s}}

\newcommand{\dimVecorsSnap}{N_x}
\newcommand{\dimParam}{p}

\newcommand{\SnapMatrix}{\bm{S}}

\newcommand{\trainedSnapMat}[1]{\bm{S}_{_{\trainingparam{#1}}}}
\newcommand{\InterpSnapMat}{\bm{S}_{_{\untrainedparam}}}

\newcommand{\mcal}[1]{\mathcal{#1}}
\newcommand{\mbb}[1]{\mathbb{#1}}

\newcommand{\txt}[1]{\textnormal{#1}}

\newcommand{\matrixSpace}[1]{\mathbb{R}^{N_{#1}\times \nbrModes}}

\newcommand{\StiefelmanifoldSpace}[1]{\mathcal{ST} (\nbrModes , N_{#1})}
\newcommand{\GrassmanifoldSpace}[1]{\mathcal{G}(\nbrModes, N_{#1})}

\newcommand{\TangsubSpaceGrass}[2]{\mcal{T}_{_{[{\Mat{#1}}_{#2}]}}\GrassmanifoldSpace{}}

\newcommand{\curve}{\Gamma}

\newcommand{\OrthGroupe}[1]{\bm{\mcal{O}}(#1)}

\newcommand{\Exp}[1]{\bm{\txt{Exp}}_{_{[#1]}}}
\newcommand{\Log}[1]{\bm{\txt{Log}}_{_{[#1]}}}

\newcommand{\Mat}[1]{{#1}}
\newcommand{\MatCalib}[1]{{#1}}

\newcommand{\distGrass}[2]{\txt{\textit{dist}}_{\mcal{G}} (#1,#2)}

\newcommand{\statevarP}{p}

\newcommand{\trainingparam}[1]{\ctrlvar_{#1}}
\newcommand{\untrainedparam}{\widetilde{\ctrlvar}}

\newcommand{\InitVel}[1]{\mcal{X}_{_{#1}}}

\newcommand{\InterpGrassVel}{\mcal{X}_{_{\untrainedparam}}}

\newcommand{\CalibMatrix}[2]{\Mat{#1}_{_{#2}}}
\newcommand{\nonintrusivemethod}{\textit{Bi-CITSGM}}

\newcommand{\statefluctvecU}{\statevecU '}
\newcommand{\statefluctvarP}{p'}

\usepackage[a4paper, total={5.5in, 9.in}]{geometry}

\usepackage[nomarkers]{endfloat} % Tables and figures at the end of manuscript

\SetKwInput{KwOffline}{Offline}
\SetKwInput{KwOnline}{Online}
\usepackage[export]{adjustbox}
\begin{document}

\begin{frontmatter}


% \author{Name\corref{cor1}\fnref{label2}}
% \address{Address\fnref{label3}}

\author{M. Oulghelou\fnref{label1}}
%\ead{}
%\tnotetext[label1]{}
\fntext[label1]{mourad.oulghelou@univ-lr.fr}
%\cortext[cor1]{}
\author{C. Allery\fnref{label2}}
%\tnotetext[label2]{}
\fntext[label2]{cyrille.allery@univ-lr.fr}
%\cortext[cor2]{}

\address{LaSIE, UMR-7356-CNRS, Universit\'e de La Rochelle P\^ole Science et Technologie,
Avenue Michel Cr\'epeau, 17042 La Rochelle Cedex 1, France.}

\title{Non intrusive method for parametric model order reduction using a bi-calibrated interpolation on the Grassmann manifold}


%\author{Mourad Oulghelou\corref{cor1}, Cyrille Allery\corref{cor2}}
% \address[a]{mourad.oulghelou@univ-lr.fr}
% \address[b]{cyrille.allery@univ-lr.fr}
% \address{}

\begin{abstract}
Approximating solutions of non-linear parametrized physical problems by interpolation presents a major challenge in terms of accuracy.
In fact, pointwise interpolation of such solutions is rarely efficient and leads generally to incorrect results. 
However, instead of using a straight forward interpolation on solutions, reduced order models (ROMs) can be interpolated.
%instead of using straight forward interpolation of solutions, 
%To overcome this issue, it is possible to interpolate reduced order models (ROMs). 
More particularly, Amsallem et al. \cite{Amsallem} proposed an efficient POD (Proper Orthogonal Decomposition) reduced order models interpolation technique based on differential geometry tools. This approach, named in this paper ITSGM (Interpolation On a Tangent Space of the Grassmann Manifold), allows through the passage to the tangent space of the Grassmann manifold, to approximate accurately the reduced order basis associated to a new untrained parameter. This basis is used afterwards to build the interpolated ROM describing the temporal dynamics by performing the Galerkin projection on the high fidelity model. 
Such Galerkin ROMs require to access to the underlying high fidelity model, leading by that to intrusive ROMs.
In this paper, and contrary to the ITSGM/Galerkin approach, we propose a non-intrusive reduced order modeling method which is independent of the governing equations. This method is named through this paper Bi-CITSGM (Bi-Calibrated ITSGM).
%This new method combines ITSGM with an optimization problem serving for modes classification. 
It consists first to interpolate the spatial and temporal POD sampling bases considered as representatives of points on Grassmann manifolds, by the ITSGM method. Then, the resulting bases modes are reclassified by introducing two orthogonal matrices. These calibration matrices are determined as analytical solutions of two optimization problems.
Results on the flow problem past a circular cylinder where the parameter of interpolation is the Reynolds number, show that for new untrained Reynolds number values, the developed approach produces satisfyingly accurate solutions in a real computational time.

\end{abstract}

\begin{keyword}
Non-intrusive Reduced Order Models (ROMs), Proper Orthogonal Decomposition (POD), bases interpolation, Grassmann manifold.
\end{keyword}

\end{frontmatter}
\newpage
% %%%%%%%%%%%%%%%%%%%%%%%%%%%%%%%%%%%%%%%%%%%%%%%%%%%%%%%%%%%%%%%%%%%%%%%%%%%%%%%%%%%%%%%%%%%%%%%%%%%%%%%%%%%%%%%%%%%%%%%%%%%%%%%%%%%%%%%%%%%%%%%%%%%
\section{Introduction}
Parametrized physical problems arise in many practical engineering design and analysis problems such as aerospace engineering \cite{Athan2004, Bloor1995}, biomedical engineering \cite{Mili2004}, flow control \cite{LU2017} etc.
The numerical computational cost of these applications can be exorbitant. 
Particularly, for applications requiring multiple resolutions of the high fidelity model over the parameter space, using a fine spatial mesh and small time-steps.
%when the solution of the underlying parametrized high fidelity model is required to be accurately solved for multiple parameters using a fine spatial mesh and small time-steps. 
For complex non-linear physical phenomena, even with modern computers, these numerical computations still remain time-consuming. Beside numerical computations, experiments can be carried out. Nevertheless, these last remain expensive  as well.
%
%with a very reduced computational budget
In order to computationally tackle the challenging cost of this kind of problems, alternative approaches must be used. 
%The challenge is to reliably achieve a very good predictability of the problem solution with a reduced computational budget
%
%These challenging costs motivate the quest for alternative ways to computationally tackle the problem. 
%This can be achieved by representing the solutions in a much lower dimensional subspace that retain the mean features contained in the original system. 
A possible approach can be by representing the solutions in a much lower dimensional subspace that retain the mean features contained in the original system. 
%
%Thus, there is a need for efficient numerical methods enabling to approximate the solutions on reduced dimensional subspaces described through a few basis functions. 
This requires efficient numerical methods that allow to compute basis functions that represent these reduced subspaces.
Such techniques, called reduced order model (ROM) methods, lead to significant computational cost savings. Especially in applications requiring multiple evaluations of the solution over the parameter space \cite{Bergmann, Tallet-Allery-2016, Tallet-Allery-2015}.
The Proper Orthogonal decomposition (POD) method is the most used method for reducing problems dimensionality. Starting from a set of available snapshot solutions at different time instants, the POD method aims to construct a subspace (represented by a spatial basis) where each snapshot solution can accurately be represented. The most striking property of POD is its optimality \cite{Berkooz-1993}. More specifically, the POD basis can be formed by eliminating the higher order POD modes carrying small energy contributions. This provides a powerful means for constructing the lowest possible dimensional basis carrying the quasi-totality of information contained in the original set of snapshots. The temporal dynamics is afterwards determined through the Galerkin projection of the high fidelity model onto the spatial basis \cite{ROWLEY2004115}.
%
%Galerkin projection \cite{ROWLEY2004115} is the most frequently used approach for obtaining a ROM from a POD basis for unsteady problems. 
%
%Essentially, the governing time-dependent equations are projected onto the POD basis forming by that
%
The resulting equations form a system of low order ordinary differential equations. It's worth noting that a POD basis is optimal only for the solution snapshots considered for its construction. Consequently, a POD basis constructed for a given parameter and used to predict the ROM for another parameter, often results in solutions lacking accuracy.
%
%Consequently the ROM built in this way may lack accuracy \cite{OULGHELOUAMC2018}, especially when a POD basis constructed for a given parameter is used to predict the ROM for another parameter. 
%
To overcome this issue, it is possible to use the ITSGM (Interpolation on Tangent spaces of the Grassmann Manifold) method based on tools of differential geometry. Let's consider an ensemble of distinct trained parameters for which a set of POD spatial bases of the same dimension is available. These bases can be seen as representatives of points on the Grassmann manifold in which geodesic paths can be associated to second order ordinary differential equations. Thus, the geodesic paths between an arbitrary chosen reference point and the other points differ from each other only from their initial velocities lying in the tangent space of the Grassmann manifold at the reference point.   
 The key idea of the ITSGM method is to interpolate the initial velocity for a new untrained parameter and use it to determine the interpolated spatial basis given as the final point of the associated geodesic path. 
%
%the ITSGM uses calculus of geodesics on the Grassmann manifold predict a spatial basis associated to a new untrained parameter. 
Finally, the interpolated ROM is built through the Galerkin projection of this basis onto the high fidelity model.
This powerful method was firstly introduced by Amsallem et al. \cite{Amsallem} in the context of aeroelasticity and applied successfully to sub-optimal control \cite{Oulghelou2017, OULGHELOUAMC2018}.
%Eventhough its effectiveness, this approach suffers from intrusiveness, which make it inapplicable for experimental data. 
%
Eventhough its effectiveness, this approach and others based on Galerkin projection suffer from intrusiveness, which make them inapplicable for experimental data. 
%Moreover, computing the POD projection coefficients requires the evaluation of integrals over the entire spatial domain which can be computationally expensive for large problems. 
%This provided a key motivation for the development of a new non intrusive version of the ITSGM method.
This provided a key motivation for the development of alternative model reduction techniques called non intrusive reduced order models.
\\
To construct a parametrized nonintrusive ROM, Xiao et al. \cite{Xiao2016,Xiao2017} suggest to use a two-level RBF interpolation. In the first RBF level, and for an arbitrary untrained parameter, the interpolated set of snapshots and their associated POD basis are calculated. In the second RBF level, the temporal dynamics is represented as a set of hyper surfaces approximated by the RBF interpolation. 
% the building of a set of hyper-surfaces that represent temporal dynamics, these hyper-surfaces are afterwards interpolated on the Smolyak grids using the RBF (Radial basis functions) method, whereas the reduced spatial basis is obtained by applying the POD to the snapshots also obtained after the RBF interpolation. 
The approach was used to study the two dimensional flow around a cylinder when the Reynolds number varies. Still on the same study case, Shinde et al. \cite{SHINDE2016} propose to simply approximate the spatial and temporal bases functions by linear interpolation of their modes. This approach was also applied by Joyner \cite{Joyner2004} for eddy current damage detection. 
%In this case, it was observed that the results obtained by the POD/interpolation method are comparable to those obtained by the POD/Galerkin method.
%
%
%
In this article, we propose a novel non-intrusive model reduction method addressed to approximate time-dependent parametrized non-linear physical phenomena. Being inspired from ITSGM, the new model reduction method is called throughout the paper Bi-CITSGM (Bi-Calibrated Interpolation on Tangent spaces of the Grassmann Manifold). In contrast to Galerkin projection based version of ITSGM, which requires the access to the underlying mathematical model, the Bi-CITSGM is non intrusive. That means that it only needs to explore the information content of the existing parametrized data snapshots to be able to approximate the solution snapshots  for a new untrained parameter. This enables the straightforward application of this methodology to generalized parametrized non-linear physical phenomena, either provided by numerical computations or experimental observations.
Starting from an ensemble of solution snapshots obtained from experiment or from solving the high fidelity model at a finite set of trained points in the parameter space, the Bi-CITSGM is carried out in two phases. The offline phase of data processing in which each ensemble of snapshots is represented in a low dimensional subspace through the POD method. The online phase of the bi- calibrated interpolation consisting first of interpolating singular values by using spline cubic; then interpolating the spatial and temporal bases by using ITSGM method; and finally, introducing according to a well defined criterion two small sized matrices serving for modes calibration. These undetermined calibration matrices obey the orthonormality condition and are determined as the solution of two separate optimization problems that aim to capture the POD modes character from the sampling trained POD bases. 
%
%The Bi-CITSGM is computationally efficient. Its computational complexity is linear function of the the number of degrees of freedom of the underlying high fidelity model. It also pointed out that in the case of affinely dependent parametric reduced order model, this complexity can significantly reduce to be proportional to the cube of the number of reduced variables. 
%
%is shown that its computational complexity can be made proportional to the cube of the number of reduced variables.
%
%it has been shown by ??? \cite{} that the complexity of the ITSGM can be made proportional to the cube of the number of reduced variables. 
%
%In terms of computational complexity, it has been shown by ??? \cite{} that the complexity of the ITSGM can be made proportional to the cube of the number of reduced variables. 
%
%This gives rise to an interpolatory model reduction technique that
%
%
%
%The new approach is addressed also the kind of problems when data are provided by experiment 
%
%The predictability of Bi-CITSGM depends on the initial choice of training solution snapshots.
%
%It has been shown by ??? \cite{} that the complexity of the ITSGM can be made proportional to the number of reduced variables. 
\\
%The remainder of this article is organized as follows: First, a brief overview of the POD method is given in section 2. 
%In Section 3, we outline the problem of bases interpolation and describe the ITSGM method. 
%In the section that follows, and under certain conditions, we show how the new non intrusive approach Bi-CITSGM is constructed. 
%A numerical example is considered in Section 4. 
%Finally, the conclusion
%is given in Section 5.
%
%
%
%
%
%
The remainder of this article is organized as follows: First, a brief overview of Geodesic paths parametrization on the Grassmann manifold is given in section 2. 
In section 3, we outline the problem of bases interpolation and present the ITSGM method for bases interpolation.
In section 4, the new proposed non intrusive approach Bi-CITSGM is formalized. 
Thereafter, the numerical example of the parametrized flow past a circular cylinder is considered in Section 5. 
Finally, the conclusion is given in Section 6.
\section{Some preliminary results of the Grassmann manifold}
In this article we are concerned with parametrized points on the Grassmann manifold \cite{Boumal-2015,Connie-2016}. In recent years, the Grassmann manifold has attracted great interest in various applications, such as subspace tracking \cite{Anuj-2004}, sparse coding \cite{Harandi-2011, Harandi2016}, clustering \cite{Cetingul-2009}, and model reduction \cite{Amsallem, OULGHELOUAMC2018}. The Grassmann manifold $\GrassmanifoldSpace{}$ is defined as the set of all $\nbrModes$-dimensional subspaces in $\mbb{R}^{N}$, where $0 \leq \nbrModes \leq N$. It is a matrix manifold that is locally similar to an Euclidean space around each of its points. A concrete representation of the Grassmann manifold is the Stiefel Manifold representation \cite{Edelman1998} given as follows
$$\GrassmanifoldSpace{} \cong \StiefelmanifoldSpace{} / \OrthGroupe{\nbrModes} $$
where $\OrthGroupe{\nbrModes}$ is the group of all $\nbrModes\times\nbrModes$ orthogonal matrices and $\StiefelmanifoldSpace{} = \{ \Xbasis{} \in \matrixSpace{} : \Xbasis{}^T \Xbasis{} = \Mat{I}_{\nbrModes} \}$ is the set of all bases of dimension $\nbrModes$, called Stiefel manifold \cite{Edelman1998, Absil}.
A point  $[\Xbasis{}] \in \GrassmanifoldSpace{}$ is defined by
$$ [\Xbasis{}] = \{ \Xbasis{} \MatCalib{Q} \mid \Xbasis{}^T\Xbasis{} = \Mat{I}_{\nbrModes}, \ \  \MatCalib{Q}\in\OrthGroupe{\nbrModes}\}$$ 
where $\Xbasis{}$ is a point on the Stiefel manifold, i.e., a matrix with orthogonal columns. $[\Xbasis{}]$ is then realized as a representative $\Xbasis{}$ from the equivalent class $\{ \Xbasis{} \MatCalib{Q} : \txt{ for all } \MatCalib{Q}\in\OrthGroupe{\nbrModes}  \}$.
%\subsection{Tangent space and geodesic distance}
%by $\Xbasis{}$ the $N$ by $\nbrModes$ full rank column matrix associated to a parameter $\trainingparam{} \in \mbb{R}^{\dimParam}$, $\dimParam \geq 1$, and
\\
At each point $[\Xbasis{}]$ of the manifold $\GrassmanifoldSpace{}$, a tangent space \cite{Edelman1998, Absil} of the same dimension \cite{Edelman1998} can be defined. 
The tangent space at a point $[\Xbasis{}]\in\GrassmanifoldSpace{}$ is given by the following  abstract concrete set \cite{Boumal-2015}
%We denote the tangent space at a point $[\Xbasis{}]\in\GrassmanifoldSpace{}$ by $\TangsubSpaceGrass{\Xbasis{}}{}$. The tangent space is given by the following  abstract concrete set \cite{Boumal-2015} 
$$ \TangsubSpaceGrass{\Xbasis{}}{} = \{ \InitVel{} \in \matrixSpace{} \mid \ \ \Xbasis{}^T \InitVel{} = 0 \}.$$
%
%
%The geodesic distance between two points $[\Mat{U}_1],[\Mat{U}_2] \in \GrassmanifoldSpace{}$ is the distance that minimises the geodesic path
%
The geodesic distance $\distGrass{\Xbasis{}}{\XbasisEnd{}}$ between two subspaces $[\Xbasis{}]$ and $[\XbasisEnd{}]$ on the Grassmann manifold is the minimum of the lengths of paths between them. By taking the SVD (Singular Value Decomposition) of $\Xbasis{}^T \XbasisEnd{}$ as $\Mat{U} \Sigma \Mat{V}^T = \Xbasis{}^T \XbasisEnd{}$ such that $\Sigma = \txt{diag}(\sigma_i)$, the geodesic distance between $[\Xbasis{}]$ and $[\XbasisEnd{}]$ is defined as the summation of squared principal angles
\begin{equation}
\distGrass{\Xbasis{}}{\XbasisEnd{}} = \sqrt{\displaystyle{\overset{}{\underset{i}{\sum}}} \arccos^2(\sigma_i)} 
\end{equation} 
\\
A path that minimizes this distance is called geodesic \cite{Wald}. This path is associated with a second order differential equation \cite{Absil, Edelman1998, Boothby} uniquely defined by two initial conditions. These are typically the initial position and initial velocity  \cite{Wald}. A geodesic path can then be represented by a twice differentiable function $\curve : [0,1] \longrightarrow \GrassmanifoldSpace{}$, where $\curve(0)$ and $\curve(1)$ are respectively its initial and final points. The parametric representation of the geodesic path in $\GrassmanifoldSpace{}$ with initial conditions $\curve(0) = [\Xbasis{}]$ and $\deriveedot{\curve}(0) = \InitVel{}
 $ is given by \cite{Absil, Edelman1998}
\begin{equation}\label{Eq:Geodesic_path}
 \curve(t) = \txt{span} \left\{ \Xbasis{} V \cos(t \Sigma ) + U \sin(t \Sigma )  \right\} \ \ \ \ 0\leq t \leq 1
\end{equation}
where $U \Sigma V^T$ is the thin SVD of the initial velocity $\InitVel{}$. For every point $[\Xbasis{}]$ in $\GrassmanifoldSpace{}$ there exist a unique geodesic starting from $[\Xbasis{}]$ in every direction $\InitVel{} \in \StiefelmanifoldSpace{}$, giving us the exponential map $\Exp{\Xbasis{}} : \TangsubSpaceGrass{\Xbasis{}}{} \longrightarrow \GrassmanifoldSpace{}$ explicitly calculated as $\Exp{\Xbasis{}}(\InitVel{}) = \curve(1) = [\XbasisEnd{}] $. The exponential of $\InitVel{}$ is given by
\begin{equation}\label{Eq:Exp_map}
[\XbasisEnd{}] = \txt{span} \{ \Xbasis{} V \cos(\Sigma) + U \sin(\Sigma)  \}
\end{equation}
Let us denote $\Log{\Xbasis{}}$ the inverse map to $\Exp{\Xbasis{}}$, which is defined only in a certain neighbourhood of $[\Xbasis{}]$. If $\Exp{\Xbasis{}}(\InitVel{}) = [\XbasisEnd{}]$, then $\InitVel{}$ is the vector determined as follows
\begin{equation}\label{Eq:Log_map}
\InitVel{} = \Log{\Xbasis{}}([\XbasisEnd{}]) = U \arctan(\Sigma) V^T
\end{equation}
where $U\Sigma V^T$ is the thin SVD of $(I-\Xbasis{} \Xbasis{}^T)\XbasisEnd{}(\Xbasis{}^T\XbasisEnd{})^{- 1}$ and $\Log{\Xbasis{}}([\Xbasis{}]) = 0$
%\subsection{Exp and Log mapping of the Grassmann manifold}
\section{Subspaces interpolation using ITSGM method}
Let $\left\{ \ctrlvar_i \in \mbb{R}^{\dimParam}, i=1,\cdots,\nbrParam\right\}$ be a set of parameters and $[\Sbase{1}], [\Sbase{2}], \dots, [\Sbase{\nbrParam}]$ the corresponding set of $\nbrParam$ parametrized subspaces \footnote{Here $\Sbase{i}$ is refereed to as a possible  representative basis of the subspace $\Sbase{i} \in \GrassmanifoldSpace{}$. Typically, this basis can be a POD basis associated to the parameter $\trainingparam{i}$.} belonging to $\GrassmanifoldSpace{}$ . Consider the problem of interpolation in which we seek an approximation of $[\InterpSbase]$ for an untrained parameter $\untrainedparam \notin \{ \trainingparam{1},\dots ,\trainingparam{\nbrParam}\}$. 
Obviously, usual processes of interpolation fail in this particular case. 
The reason is that the Grassmann manifold is not a flat space. Hence, a straight forward interpolation of its points does not necessarily result in a point that is included in it. It is therefore necessary to reformulate the interpolation process to be suitable for points in the Grassmann manifold. A possible way to do it is by performing interpolations on initial velocities lying on a tangent space (flat space) at a point of $\GrassmanifoldSpace{}$. The ITSGM (Interpolation on a Tangent space of the Grassmann Manifold) procedure proposed in \cite{Amsallem} helps to complete this task. The steps of ITSGM are given as follows
\begin{itemize}
\item[\hspace{10pt}]
\begin{itemize}
\item[ \textit{step 1}] Choose the origin point of tangency, for example $[\Sbase{i_0}]$ where $i_0 \in \{1,\dots,\nbrParam\}$.
\item[ \textit{step 2}] For $i \in \{1,\dots,\nbrParam\}$, map the point $[\Sbase{i}] \in \GrassmanifoldSpace{}$ to $\InitVel{i} \in \TangsubSpaceGrass{\Xbasis{i_0}}{}$ such that $\InitVel{i} = \Log{\Sbase{i_0}}(\Sbase{i})$ is the vector represented by 
$$ \InitVel{i} = U_i \arctan(\Sigma_i) V_i^T $$
where $U_i \Sigma_i V_i^T = (I-\Sbase{i_0} \Sbase{i_0}^T)\Sbase{i}(\Sbase{i_0}^T\Sbase{i})^{- 1}, i=1, \dots, \nbrParam$, are thin SVD.
\item[ \textit{step 3}] Interpolate the initial velocities $\InitVel{1}, \InitVel{2}, \dots, \InitVel{\nbrParam}$ for the untrained parameter $\untrainedparam$ using a standard interpolation and obtain $\InterpGrassVel$.
\item[ \textit{step 4}] Finally by the exponential mapping, map the interpolated velocity $\InterpGrassVel$ back to the Grassmann manifold. The matrix representation of the interpolated subspace is given by
$$
    \InterpSbase = \Sbase{i_0} \widetilde{V} \cos(\widetilde{\Sigma}) + \widetilde{U} \sin(\widetilde{\Sigma})
$$
where $\widetilde{U} \widetilde{\Sigma} \widetilde{V}^T$ is the thin SVD of the initial velocity vector $\InterpGrassVel$.
\end{itemize}

\end{itemize}
For more details the reader can refer to Amsallem et al. \cite{Amsallem} and Absil et al. \cite{Absil}.

\section{Bi-CITSGM method}
	\subsection{Problem statement}
		Given a set of parameters $\left\{ \trainingparam{i} \in \mbb{R}^{\dimParam}, i=1,\cdots,\nbrParam\right\}$ for which we dispose of $\nbrParam$ solutions at $\nbrSnap$ different time instants of a discretized time dependant non-linear parametrized physical problem of size $\dimVecorsSnap$. For each parameter, these solutions can be stored in a $\dimVecorsSnap\times\nbrSnap$ matrix.
The aim is to approximate the solution snapshots matrix $\InterpSnapMat$ for a new untrained parameter $\untrainedparam \notin \{ \trainingparam{1},\dots,\trainingparam{\nbrParam}\}$ without going through costly numerical or experimental processes.
A trivial approach to solve this problem is to use standard interpolation approaches. 
These methods such as Lagrange or RBF (Radial Basis Function) are known to be efficient for linear or weakly non-linear problems, while for non-linear phenomena, they are pointless and unuseful.
The aim of this paper is thus to overcome this limitation by introducing a new interpolation approach based on the ITSGM, which takes into account the non-linear behaviour of the physical problem. In the following, this interpolation method is 
%is based on differential geometry tools and 
called Bi-CITSGM (Bi-Calibrated Interpolation on Tangent spaces of the Grassmann manifold).

	\subsection{Description of the method}		
		Let us assume that $\trainedSnapMat{i}$ is approximated in the low dimensional subspace $[\Sbase{i}]$ calculated by the POD method \footnote{The POD basis can be chosen to be optimal with respect to any inner product, namely, $L^2$ or $H^1$ or Euclidean inner product. For the sake of simplicity, the Euclidean 2-norm is used in the following.} as follows
\begin{equation}\label{POD approx sampling}
\trainedSnapMat{i} \approx \Sbase{i} \trainedSingVal{i} \Tbase{i}^T, \hspace*{0.3cm} i=1,\dots,\nbrParam,
\end{equation}
where $\trainedSingVal{i}\in \mbb{R}^{\nbrModes\times \nbrModes}$ is the matrix of singular values (square roots of POD eigenvalues) of $\trainedSnapMat{i}$, and $\Sbase{i}\in \StiefelmanifoldSpace{x}$ and $\Tbase{i} \in \StiefelmanifoldSpace{s}$, $ i=1,\dots,\nbrParam$, are the corresponding spatial and temporal bases. 
The first step of the $\nonintrusivemethod$ method is to approximate the matrix of singular values $\InterpSingVal$ by spline cubic interpolation. 
%This can be done easily trough a usual interpolation approach, such as Spline, RBF, or Lagrange interpolations. 
Next, by using the ITSGM method, the subspaces $[\Sbase{i}]$ for $i=1,\dots,\nbrParam$ (resp. $[\Tbase{i}]$) are interpolated to obtain the approximated subspace $[\InterpSbase]$ (resp. $[\InterpTbase]$). Since $\InterpSbase$ and $\InterpTbase$ are only two possible representative bases of the subspaces $[\InterpSbase]$ and $[\InterpSbase]$, the right classification of modes according to the interpolated POD eigenvalues is not verified.
Thus, the approximated snapshots matrix $\InterpSnapMat$ must be calibrated as follows 
\begin{equation*} 
\InterpSnapMat = \InterpSbase \CalibMatrix{Q}{x} \InterpSingVal \CalibMatrix{Q}{t}^T \InterpTbase^T
\end{equation*}
where $\CalibMatrix{Q}{x},\CalibMatrix{Q}{t} \in \left[\mcal{O}(q)\right]^2$ are introduced to recover the modes POD character of $\InterpSbase$ and $\InterpTbase$ respectively. These two matrices are of a big significance in the expression of $\InterpSnapMat$. In fact, their introduction guarantees the adequate combination of column vectors of $\InterpSbase$ and $\InterpTbase$ which realizes the best fit with the interpolated singular values matrix $\InterpSingVal$.
\\
The question of finding the best calibration matrices $\CalibMatrix{Q}{x}$ and $\CalibMatrix{Q}{t}$ is a delicate one. 
%This paper is dedicated to propose a practical approach to solve this problem. 
First of all, the modes signs of sampling spatial and temporal bases must be adjusted with respect to a reference basis \footnote{In this manipulation, signs of spatial and temporal modes are simultaneously adjusted. For example, the $j^{\txt{th}}$ modes associated to $\Sbase{k}$ and $\Tbase{k}$ can be multiplied by $-1$, without causing any effect on the snapshot matrices representation given by \eqref{POD approx sampling}.}. The reference basis denoted $\Sbase{k_0}$, is chosen as the closest basis to $\InterpSbase$ in the sense of the Grassmannian distance. Thus, $k_0$ is determined by
$$ k_0 = \underset{i\in\{1,\dots,\nbrParam\}}{\txt{\textbf{argmin}}} \ \ \distGrass{\InterpSbase}{\Sbase{i}} $$
Denote by $||\cdot||_2$ the Euclidean two-norm. Once $\Sbase{k_0}$ is known, the signs adjustment is performed by algorithm \ref{Alg:ajust_signs_Bi-CITSGM}.
\\
\begin{algorithm}[H]
 \SetAlgoLined
\For{$k=1,\dots,\nbrParam$}{
\For{$j=1,\dots,\nbrModes$ \textbf{and} $j\neq j_0$}{
	{\small If $ || \Sbase{k_0}^j - \Sbase{k}^j ||_2 > || \Sbase{k_0}^j + \Sbase{k}^j ||_2$ \\
	\hspace*{0.4cm} multiply the $j^{\txt{th}}$ spatial and temporal modes $\Sbase{k}^j$ and $\Tbase{k}^j$ by $-1$.}	
	}	
}
\caption{Bases modes signs adjustment.}
\label{Alg:ajust_signs_Bi-CITSGM}
\end{algorithm}
In the proposed non intrusive model reduction approach Bi-CITSGM, the steps of determining $\CalibMatrix{Q}{x}$ and $\CalibMatrix{Q}{t}$ are the same. Thus, the methodology is described only for the finding of the spatial calibration matrix $\CalibMatrix{Q}{x}$. 
The same steps applies to determine the temporal calibration matrix $\CalibMatrix{Q}{t}$.
%Now that the spatial bases are signs adjusted, we seek the calibration matrix $\CalibMatrix{Q}{x}$ determined as the solution of the constrained minimization problem
The spatial calibration matrix $\CalibMatrix{Q}{x}$ is determined as the solution of the constrained minimization problem
\begin{equation}\label{Eq : minimization_problems_Bi-CITSGM}
\underset{\CalibMatrix{Q}{x}\in \mcal{O}(q)}{\txt{min}} \ \ \somme{i}{1}{\nbrParam} \omega_i  || \InterpSbase  \CalibMatrix{Q}{x} - \Sbase{i}||_F^2
\end{equation}
where $||\cdot||_F$ is the Frobenius norm and $\omega_i$ are the spatial weights given for $m>1$ as follows
\begin{equation}\label{Bi-CITSGM_spatial_weights}
\omega_i = \myfrac{\distGrass{\InterpSbase}{\Sbase{i}}^{-m}}{\somme{k}{1}{\nbrParam} \distGrass{\InterpSbase}{\Sbase{k}}^{-m}} 
\end{equation}
%To solve the 
%
%
%
%\begin{proposition}
%Given the eigendecomposition of the matrix
%$$\somme{j}{1}{\nbrParam}  \somme{i}{1}{\nbrParam} \omega_i \omega_j\Sbase{i}^T\InterpSbase \InterpSbase^T \Mat{U}_{j} = \Mat{P} \Mat{\lambda} \Mat{P}^T $$
%The minimization problem \eqref{Eq : minimization_problems_Bi-CITSGM} admits a unique solution given by
%\begin{equation*}
%\CalibMatrix{Q}{x} = \InterpSbase^T \left( \somme{i}{1}{\nbrParam} \omega_i\Sbase{i} \right) \Mat{P} \Mat{\lambda}^{- \frac{1}{2}} \Mat{P}^T
%\end{equation*}
%\end{proposition}
%\begin{proof}
The Lagrange function associated to the constrained minimization problem \eqref{Eq : minimization_problems_Bi-CITSGM} writes
\begin{align*}
\mcal{L}(\CalibMatrix{Q}{x},\Mat{R}) 
&=  \somme{i}{1}{\nbrParam} \omega_i || \InterpSbase \CalibMatrix{Q}{x} - \Sbase{i}||_F^2 + \txt{Trace}\left[  \Mat{R}(\CalibMatrix{Q}{x}^T\CalibMatrix{Q}{x} - I_q) \right]
\end{align*}
where $\Mat{R}$ is a symmetric matrix referred to as the "Lagrange multiplier". Considering the fact that $ \omega_1+\cdots+\omega_{\nbrParam} = 1$, the Lagrange function simplifies to
\begin{align*}
\mcal{L}(\CalibMatrix{Q}{x},\Mat{R}) 
&=  \somme{i}{1}{\nbrParam} \omega_i || \InterpSbase \CalibMatrix{Q}{x}||_F^2  - 2 \somme{i}{1}{\nbrParam} \omega_i \txt{Trace}\left( \Sbase{i}^T \InterpSbase \CalibMatrix{Q}{x} \right) +  \somme{i}{1}{\nbrParam} \omega_i || \Sbase{i}||_F^2 + \txt{Trace}\left[  \Mat{R}(\CalibMatrix{Q}{x}^T\CalibMatrix{Q}{x} - I_q) \right]
%\\
%&=  
%\txt{Trace}\left( \Mat{ K}^T \underbrace{\InterpSbase^T \InterpSbase}_{I_{\nbrModes}} \Mat{ K} \right)  - 2 \somme{i}{1}{\nbrParam} \omega_i \txt{Trace}\left( \Sbase{i}^T \InterpSbase \Mat{ K} \right) 
%\\
%& \hspace*{1.5cm} +  \somme{i}{1}{\nbrParam} \omega_i \txt{Trace}\left( \underbrace{\Sbase{i}^T \Sbase{i}}_{I_{\nbrModes}} \right) + \txt{Trace}\left(  \Mat{R}(\CalibMatrix{Q}{x}^T\CalibMatrix{Q}{x} - I_q) \right)
\\
&=\txt{Trace}\left[ (I_{\nbrModes} + R) \CalibMatrix{Q}{x}^T  \CalibMatrix{Q}{x} \right]  - 2 \somme{i}{1}{\nbrParam} \omega_i \txt{Trace}\left( \Sbase{i}^T \InterpSbase \CalibMatrix{Q}{x} \right) - \txt{Trace}\left(  \Mat{R} \right) +  \nbrModes
\end{align*}
Using the following differentiation identities
\begin{align*}
\myfrac{d}{d \Mat{B}} Trace(\Mat{A} \Mat{B}^T\Mat{B}) &= \Mat{B} (\Mat{A} + \Mat{A}^T) \\
\myfrac{d}{d \Mat{B}} Trace(\Mat{A} \Mat{B})&= \Mat{A}^T
\end{align*}
The differential of $\mcal{L}$ with respect to $\CalibMatrix{Q}{x}$ yields
\begin{equation}\label{Eq : differentielle Lag wrt K}
\deriveepartielle{\mcal{L}}{\CalibMatrix{Q}{x}}(\CalibMatrix{Q}{x},\Mat{R}) = 2\CalibMatrix{Q}{x}(I_{\nbrModes} + R) - 2 \InterpSbase^T \somme{i}{1}{\nbrParam} \omega_i \Sbase{i} 
\end{equation}
When $\deriveepartielle{\mcal{L}}{\CalibMatrix{Q}{x}}(\CalibMatrix{Q}{x},\Mat{R}) $ vanishes, $\CalibMatrix{Q}{x}$ is a stationary point characterized by
\begin{equation}\label{Eq: gradient_Lag_Bi-CITSGM}
\CalibMatrix{Q}{x} \left(  I_q + \Mat{R}\right)  = \InterpSbase^T \somme{i}{1}{\nbrParam} \omega_i  \Sbase{i} 
\end{equation}
By considering the SVD decomposition of the matrix 
\begin{equation}\label{Eq: SVD K(I+R)}
%\InterpSbase^T \somme{i}{1}{\nbrParam} \omega_i \Sbase{i} = \Mat{\xi} \Mat{\Theta} \Mat{\eta}^T
\CalibMatrix{Q}{x} \left(  I_q + \Mat{R}\right)  = \Mat{\xi} \Mat{\Theta} \Mat{\eta}^T
\end{equation}
It follows that
\begin{equation}
\left(  I_q + \Mat{R}\right)^T  \CalibMatrix{Q}{x}^T  \CalibMatrix{Q}{x} \left(  I_q + \Mat{R}\right)  = \Mat{\eta} \Mat{\Theta}^2 \Mat{\eta}^T
\end{equation}
Given that $\CalibMatrix{Q}{x}$ is orthonormal and that $\Mat{R}$ is symmetric, we can write
\begin{equation}
I_q + \Mat{R}  = \Mat{\eta} \Mat{\Theta} \Mat{\eta}^T
\end{equation}
By plugging the expression of $I_q + \Mat{R}$ into the equation \eqref{Eq: SVD K(I+R)}, the spatial calibration matrix $\CalibMatrix{Q}{x}$ is expressed as follows
\begin{equation*}
\CalibMatrix{Q}{x} =  \Mat{\xi} \Mat{\eta}^T
\end{equation*}
%\end{proof}
In the same fashion, consider the temporal weights given for $l\geq 1$ by
\begin{equation}\label{Bi-CITSGM_temporal_weights}
\kappa_i = \myfrac{\distGrass{\InterpTbase}{\Tbase{i}}^{-l}}{\somme{k}{1}{\nbrParam} \distGrass{\InterpTbase}{\Tbase{k}}^{-l}}
\end{equation}
The temporal calibration matrix $\CalibMatrix{Q}{t}$ is expressed by
\begin{equation*}
\CalibMatrix{Q}{t} =  \Mat{\zeta} \Mat{\rho}^T
\end{equation*}
where $\Mat{\zeta}$ and $\Mat{\rho}$ are the left and right singular eigenvectors of the matrix $\InterpTbase^T \somme{i}{1}{\nbrParam} \kappa_i  \Tbase{i} $.
\begin{remark}
%If the matrix $\InterpSbase^T \somme{i}{1}{\nbrParam} \omega_i \Sbase{i}$ (respectively $\InterpTbase^T \somme{i}{1}{\nbrParam} \kappa_i  \Tbase{i} $) is singular, which means that $\Mat{\Theta}$ (respectively $\delta$) contains at least one zero singular value, it
If the matrix $\Mat{A} = \InterpSbase^T \somme{i}{1}{\nbrParam} \omega_i \Sbase{i}$ is singular, we can use SVD to approximate its inverse called pseudo-inverse with the following matrix  
\begin{equation}
\Mat{A}^{-1} \approx \Mat{A}^+ = \Mat{\eta} \Mat{\Theta}^+ \Mat{\xi}^T
\end{equation}
where for a small threshold $\varepsilon>0$, the matrix $\Mat{\Theta}^+$ is given by 
\begin{equation}
\Mat{\Theta}^+ = 
\begin{cases}
1/\theta_i \hspace*{0.2cm}&\txt{if } \theta_i>\varepsilon 
\\
0 \hspace*{0.2cm}&\txt{otherwise} 
\end{cases}
\end{equation}
In this case, the spatial calibration matrix is approximated as follows
\begin{equation}
\CalibMatrix{Q}{x} \approx  \Mat{\xi} \Mat{I}_{_{\CalibMatrix{Q}{x}}}^+ \Mat{\eta}^T
\end{equation}
where 
\begin{equation}
\Mat{I}_{_{\CalibMatrix{Q}{x}}}^+ = 
\begin{cases}
1 \hspace*{0.2cm}&\txt{if } \theta_i>\varepsilon 
\\
0 \hspace*{0.2cm}&\txt{otherwise} 
\end{cases}
\end{equation}
%such that $\Mat{I}^+ = \Mat{\Theta}\Mat{\Theta}^+$.
The same remark apply to approximate the temporal calibration matrix and we have
\begin{equation}
\CalibMatrix{Q}{t} \approx  \Mat{\zeta} \Mat{I}_{_{\CalibMatrix{Q}{t}}}^+ \Mat{\gamma}^T
\end{equation}
where 
\begin{equation}
\Mat{I}_{_{\CalibMatrix{Q}{t}}}^+ = 
\begin{cases}
1 \hspace*{0.2cm}&\txt{if } \delta_i>\varepsilon 
\\
0 \hspace*{0.2cm}&\txt{otherwise} 
\end{cases}
\end{equation}
\end{remark}
The different steps of the Bi-CITSGM method are summarized in algorithm \ref{alg:Bi-CITSGM}.
\\
\begin{algorithm}[H]
\SetAlgoLined
%\underline{\textbf{\textit{Étape offline}:}}
\KwOffline{}
For each parameter $\trainingparam{i}$, $i=1,\dots,\nbrParam$, approximate the snapshots matrix $\trainedSnapMat{i}$ by the truncated SVD decomposition of order $\nbrModes$, that is
$$\trainedSnapMat{i} \approx \Sbase{i} \trainedSingVal{i} \Tbase{i}^T$$
\\
\KwOnline{}
Interpolate singular values matrices $\trainedSingVal{1}, \cdots, \trainedSingVal{\nbrParam}$ and obtain $\InterpSingVal$ the singular value matrix associated to the untrained parameter $\untrainedparam$.
\\
Interpolate $[\Sbase{i}]$ and $[\Tbase{i}]$, $i=1,\dots,\nbrParam$, by using the ITSGM method and obtain the spatial and temporal bases $\InterpSbase$ and $\InterpTbase$.
\\
Perform algorithm \ref{Alg:ajust_signs_Bi-CITSGM} to adjust the sampling bases modes signs.
\\
Calculate the weights $\omega_i$ and $\kappa_i$ using equations \eqref{Bi-CITSGM_spatial_weights} and \eqref{Bi-CITSGM_temporal_weights}.
\\
Perform SVD decompositions 
$$\InterpSbase^T \somme{i}{1}{\nbrParam} \omega_i \Sbase{i} = \Mat{\xi} \Mat{\Theta} \Mat{\eta}^T
\hspace*{1cm}
\InterpTbase^T \somme{i}{1}{\nbrParam} \kappa_i \Tbase{i} = \Mat{\zeta} \Mat{\delta} \Mat{\rho}^T$$
\\
Calculate $\Mat{I}_{_{\CalibMatrix{Q}{x}}}^+$ and $\Mat{I}_{_{\CalibMatrix{Q}{t}}}^+$ and evaluate the calibration matrices $\CalibMatrix{Q}{x}$ and $\CalibMatrix{Q}{t}$ as follows
$$\CalibMatrix{Q}{x} =  \Mat{\xi} \Mat{I}_{_{\CalibMatrix{Q}{x}}}^+ \Mat{\eta}^T
\hspace*{1cm}
\CalibMatrix{Q}{t} =  \Mat{\zeta} \Mat{I}_{_{\CalibMatrix{Q}{t}}}^+ \Mat{\rho}^T$$
\\
Reconstruction of the interpolated snapshots matrix 
$$\InterpSnapMat = \InterpSbase \CalibMatrix{Q}{x} \InterpSingVal \CalibMatrix{Q}{t}^T \InterpTbase^T$$
\caption{Algorithme Bi-CITSGM.}
\label{alg:Bi-CITSGM}
\end{algorithm} 

	\subsection{Computational complexity}	
If a univariate linear interpolation is carried out on the tangent space of the Grassmann manifold, it was demonstrated in \cite{Amsallem} that the computational cost of ITSGM performed on spatial bases is proportional to $\mcal{O}(\dimVecorsSnap \nbrModes^2)$. Similarly a computational cost proportional to $\mcal{O}(\nbrSnap \nbrModes^2)$ is expected for temporal bases interpolation. The other operations of the Bi-CITSGM algorithm involve matrix-matrix products and Singular value decompositions. Their total computational cost still proportional to $\mcal{O}(\dimVecorsSnap \nbrModes^2)$ and $ \mcal{O}(\nbrSnap \nbrModes^2)$. 
Since in general $\nbrSnap \ll \dimVecorsSnap$, the complexity of Bi-CITSGM can be approximated as being proportional to $\mcal{O}(\dimVecorsSnap \nbrModes^2)$. In other words, the complexity is a linear function of the number of degrees of freedom of the underlying higher-order computational model, which makes the Bi-CITSGM method computationally efficient.
%
%
%
%For affinely dependent parametric reduced order models, it is possible to drastically reduce the computational cost of ITSGM by considering the interpolation procedure proposed by Nguyen \cite{Nguye-2012}. In this particular case, the complexity of the Bi-CITSGM can reduce to be proportional to $\mcal{O}(\nbrModes^3)$
%
%
Still on the particular case where linear interpolation is used on the tangent space of the Grassmann manifold, it is possible to drastically reduce the computational cost of Bi-CITSGM by considering the interpolation procedure proposed by Nguyen \cite{Nguye-2012}. In this case, the complexity can reduce to be proportional to $\mcal{O}(\nbrModes^3)$.

\section{Numerical example : flow past a cylinder}
	In this section, we test the effectiveness of the Bi-CITSGM on the example of flow past a cylinder. Starting from a set of snapshots solutions calculated for different trained Reynolds number values, the goal is to interpolate these parametrized snapshots by the Bi-CITSGM, in order to obtain the solution snapshots for new untrained Reynolds number values.
\subsection{Problem description}
\begin{figure}[hbtp!]
\centering 
\includegraphics[width=0.9\linewidth]{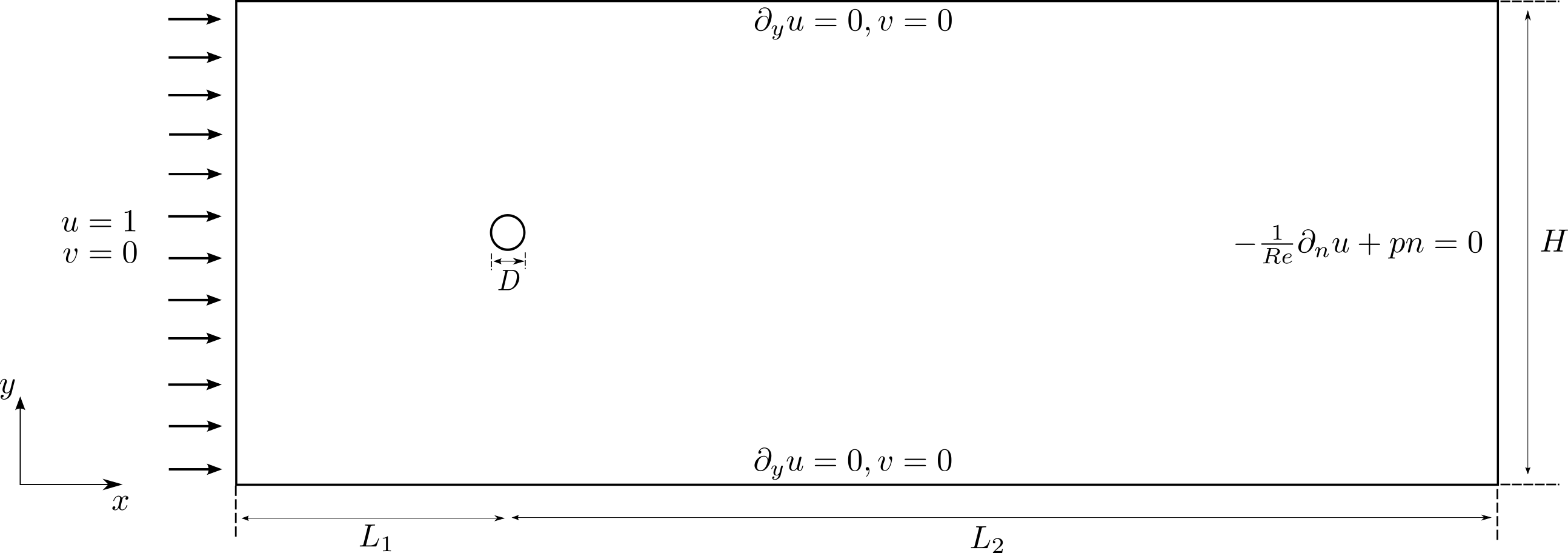}
\caption{Two-dimensional domain and boundary conditions for the problem of flow past a circular cylinder.}
\label{Fig:flow past cylinder}
\end{figure}

Consider a two dimensional flow past a circular cylinder of diameter $D$ depicted in Figure \ref{Fig:flow past cylinder}. The problem domain is rectangular with length $H = 30D$ and width  $45 D$. The center of the cylinder is situated at $L_1 = 10D$ from the left boundary and $H/2$ from the lower boundary. The fluid dynamics of the flow is driven by an inlet velocity $U$ of a unit magnitude, which enters from the left boundary of the domain, and is allowed to flow past through the right boundary of the domain. Free slip boundary conditions are applied to the horizontal edges whilst no slip boundary condition are considered on the cylinder’s wall. 
This flow problem is governed by Navier-Stokes equations given as follows 
\begin{equation}\label{EQ: Navier_Stokes_cylinder}
\begin{cases}
\deriveepartielle{\statevecU}{t} - \myfrac{1}{Re} \Delta \statevecU + \statevecU\cdot\nabla \statevecU + \nabla \statevarP  = 0  & \txt{in } \Omega\times[0,T[\\
\nabla\cdot\statevecU = 0 & \txt{in } \Omega\times[0,T[\\
\statevecU = U & \txt{on } \Gamma_{inflow} \times[0,T[\\
\statevecU = 0 & \txt{on } \Gamma_{noslip} \times[0,T[\\
\stateU_2 = 0 & \txt{on } \Gamma_{freeslip} \times[0,T[\\
\deriveepartielle{\stateU_1}{\vect{n}} = 0 & \txt{on } \Gamma_{freeslip} \times[0,T[\\
-\myfrac{1}{Re} \deriveepartielle{\statevecU}{\vect{n}} + \statevarP \vect{n} = 0  &\txt{on } \Gamma_{out}\times[0,T[\\
\statevecU(0) = \statevecU_0 & \txt{in } \Omega
\end{cases}
\end{equation}
% BEGIN FIGURE -----------------------------------------------------------------------------------------------------------------------------
% END FIGURE -------------------------------------------------------------------------------------------------------------------------------
The parameter to vary is the Reynolds number $Re = U D / \nu$ ranged from $90$ to $450$, where the jump between two training parameters is equal to $30$. 
Numerical simulations were performed with Fenics \cite{FenicsBook} (Taylor-Hood finite element $\mbb{P}_2 / \mbb{P}_1 $) using the time step $ 0.01 $ and a non-uniform mesh, including $ 85124$ DOFs for velocity and $ 10694$ DOFs for the pressure. The initial condition considered for all the training Reynolds number values corresponds to the solution at a given instant of the periodic flow regime at $ Re = 100 $. 
\subsection{Construction of the POD sampling bases}
\begin{figure}[hbtp!]
%\centering 
\hspace*{-2cm}
\includegraphics[width=1.3\linewidth]{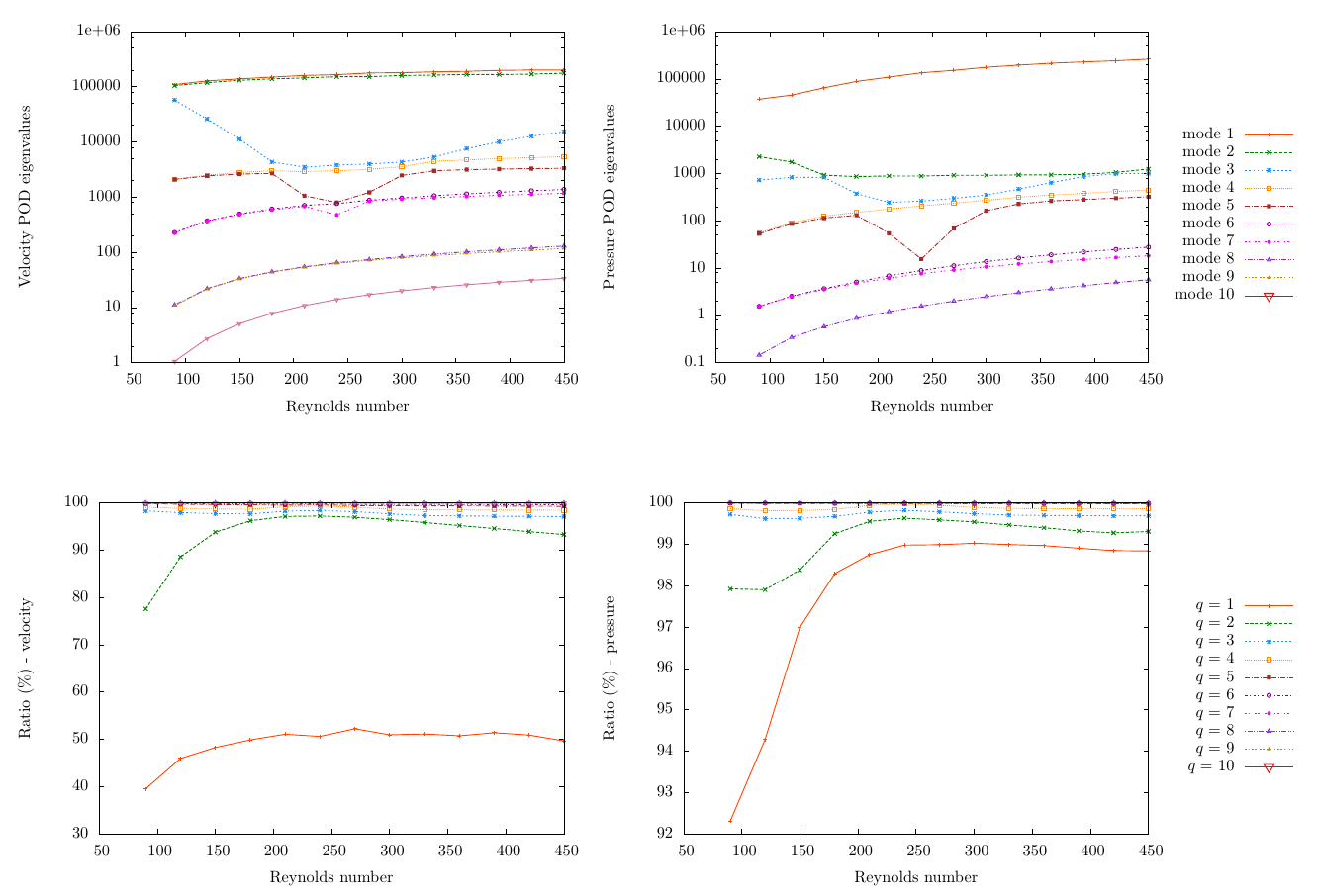}
\caption{POD eigenvalues and corresponding ratio with respect to trained Reynolds numbers. $\nbrModes$ is referred to as the number of POD modes.}
\label{Fig:POD_Eigen_values_and_Ratio_Eigenvalues}
\end{figure}

The training flows are obtained by solving the high fidelity model \footnote{In this example, the non dimensional Navier-Stokes equations are solved. The variation of the Reynolds number is carried out through the variation of the kinematic viscosity $\nu$. The inlet velocity $U$ in this case is kept constant.} for the training Reynolds number values included in $\left\{ 90, 120, 150,\dots, 420,450\right\}$ where the final time $T$ is equal to $12$. 
%
%The parameter to vary is the Reynolds number. For each Reynolds number in $ \left\{ 90, 120, 150, 180,\dots, 450\right\}$, a solution is calculated by solving equations \ref{EQ: Navier_Stokes_cylinder} in the time interval $[0,12]$. 
%
The velocity and pressure variables are decomposed into mean and fluctuating parts as follows
\begin{equation*}
\begin{cases}
\statevecU(t,x, Re) = \overline{\statevecU}(x) + \statefluctvecU(t,x, Re),  \hspace{0.5cm}\txt{ in } \Omega \\
p(t,x, Re) = \overline{p}(x) + \statefluctvarP(t,x, Re),  \hspace{0.5cm}\txt{ in } \Omega \\
%\overline{\statevecU} = U,  \hspace{0.5cm}\txt{ on } \Gamma_{inflow}\\
\end{cases}
\end{equation*}
The mean parts $\overline{\statevecU}$ and $\overline{p}$ are given by
\begin{equation*}
\overline{\vect{\stateU}} = \myfrac{1}{\nbrParam \nbrSnap} \somme{i}{1}{\nbrParam}\somme{j}{1}{\nbrSnap} \statevecU(t_j, x, Re_i)
\hspace*{2cm}
\overline{p} = \myfrac{1}{\nbrParam \nbrSnap} \somme{i}{1}{\nbrParam}\somme{j}{1}{\nbrSnap} \statevarP(t_j, x, Re_i)
\end{equation*}
where $\nbrParam$ is the number of trained Reynolds numbers and $\nbrSnap$ the number of snapshots. Two POD bases of dimensions $\nbrModes_u$ and $\nbrModes_p$ for the fluctuating velocity and pressure are then constructed and the solutions are approximated as follows
\begin{equation}\label{EQ : subspace approx POD u and p}
\statefluctvecU \approx  \somme{i}{1}{\nbrModes_u} \Tmode{u}{i} \Smode{u}{i} \hspace*{1.5cm} \statefluctvarP \approx  \somme{i}{1}{\nbrModes_p} \Tmode{p}{i} \Smode{p}{i}
\end{equation}
This two POD bases were constructed by considering $500$ snapshots ($N_s = 500$) regularly distributed
between instants $t_i = 7$ and $t_f = 12$, representing about $8$ periods of the flow. The contribution ratio of the $k^{\txt{th}}$ mode is given by
$$\txt{RIC}^{k} = \somme{i}{1}{k}\lambda_i / \somme{i}{1}{\nbrSnap}\lambda_i$$
where $\lambda_i$ are the POD eigenvalues. POD eigenvalues and respective modes contribution ratio for the different trained Reynolds number values are represented in figure \ref{Fig:POD_Eigen_values_and_Ratio_Eigenvalues}. It can  be noticed that only few modes are capable of reproducing the quasi-totality of the flow. Therefore, for this study case, the POD bases were constructed by considering the $10$ first modes for velocity and $8$ first modes for the pressure. 
In what follows, the non intrusive ROM obtained by the proposed method Bi-CITSGM is compared to the ROM obtained by ITSGM/Galerkin. 
The Galerkin based reduced order model construction for this problem is given in \ref{Appendix.Cylinder}.
% END FIGURE -------------------------------------------------------------------------------------------------------------------------------

%
\subsection{Application of the Bi-CITSGM to predict the flow for untrained Reynolds number values}
\begin{figure}[hbtp!]
\centering
\begin{subfigure}{.5\textwidth}
  \centering
  \hspace*{-4cm}\includegraphics[width=\linewidth]{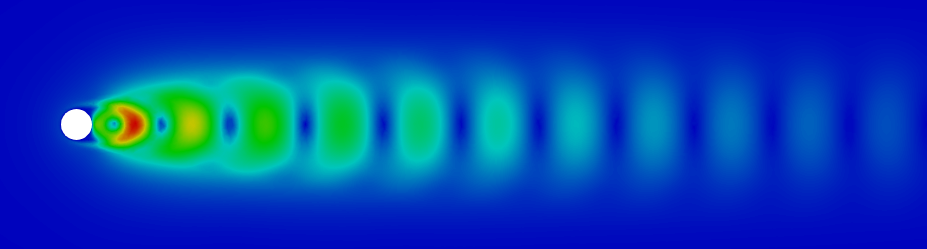}\hspace*{1cm}%
  \includegraphics[width=\linewidth]{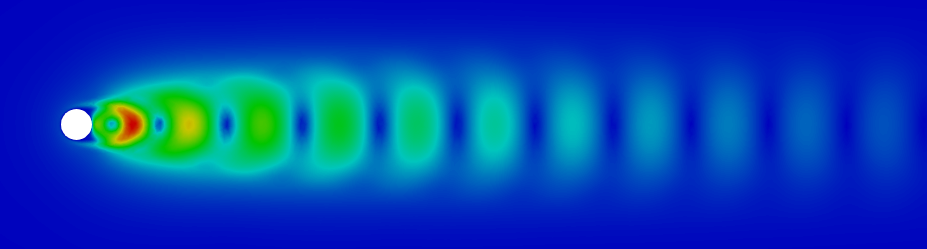}
  \caption{mode 1}
  \label{fig:sub1}
\end{subfigure}
\par\bigskip
\begin{subfigure}{.5\textwidth}
  \centering
  \hspace*{-4cm}\includegraphics[width=\linewidth]{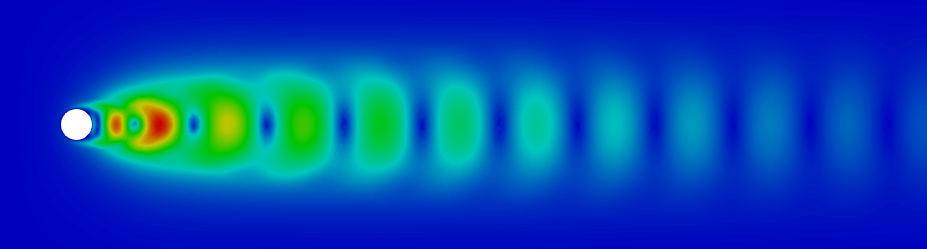}\hspace*{1cm}%
  \includegraphics[width=\linewidth]{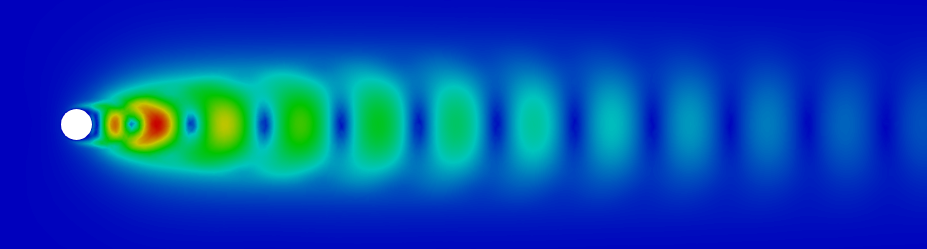}
  \caption{mode 2}
  \label{fig:sub1}
\end{subfigure}
\par\bigskip
\begin{subfigure}{.5\textwidth}
  \centering
  \hspace*{-4cm}\includegraphics[width=\linewidth]{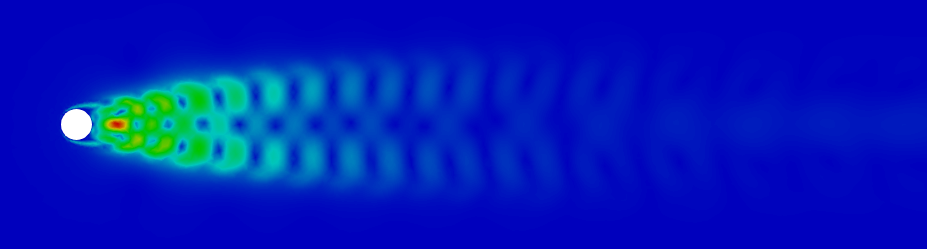}\hspace*{1cm}%
  \includegraphics[width=\linewidth]{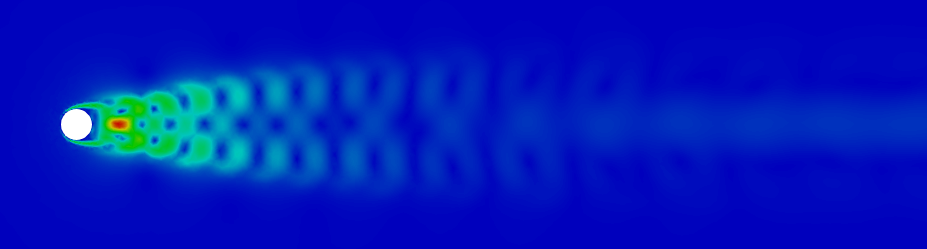}
  \caption{mode 3}
  \label{fig:sub1}
\end{subfigure}
\par\bigskip
\begin{subfigure}{.5\textwidth}
  \centering
  \hspace*{-4cm}\includegraphics[width=\linewidth]{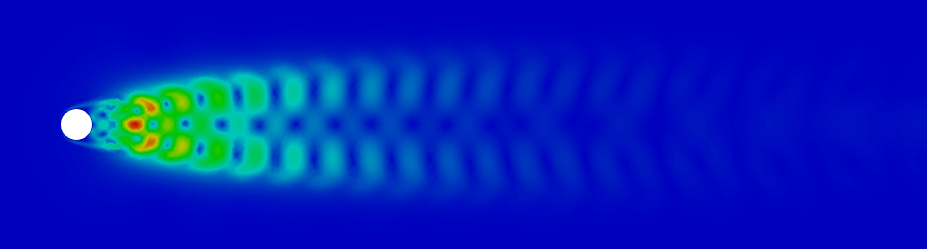}\hspace*{1cm}%
  \includegraphics[width=\linewidth]{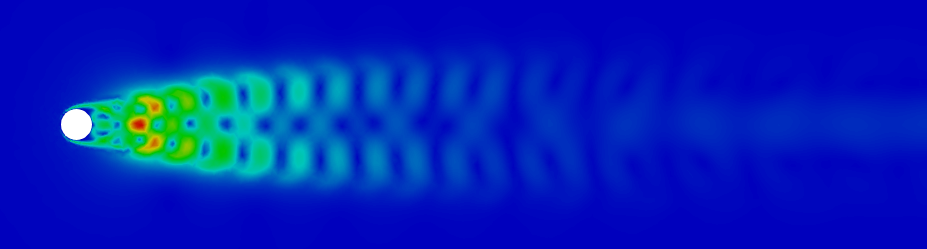}
  \caption{mode 4}
  \label{fig:sub1}
\end{subfigure}
\caption{POD (left) and calibrated ITSGM (right) velocity modes associated to the untrained Reynolds number $Re = 195$.}
\label{fig:velocity_space_modes_ITSGM_195}
\end{figure}

\begin{figure}[hbtp!]
\centering
\begin{subfigure}{.5\textwidth}
  \centering
  \hspace*{-4cm}\includegraphics[width=\linewidth]{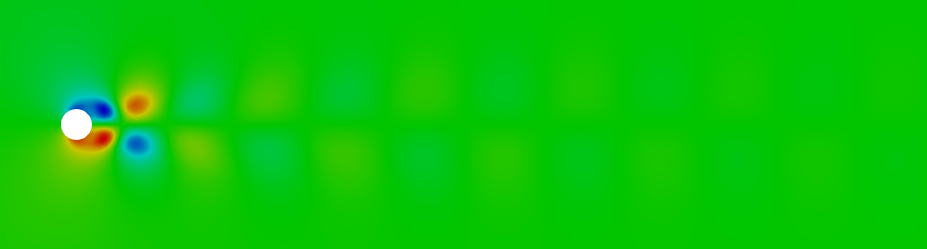}\hspace*{1cm}%
  \includegraphics[width=\linewidth]{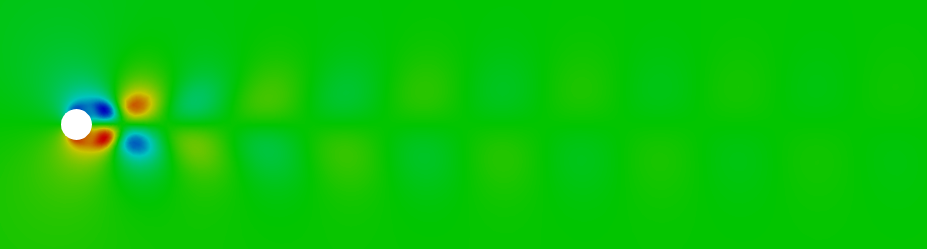}
  \caption{mode 1}
  \label{fig:sub1}
\end{subfigure}
\par\bigskip
\begin{subfigure}{.5\textwidth}
  \centering
  \hspace*{-4cm}\includegraphics[width=\linewidth]{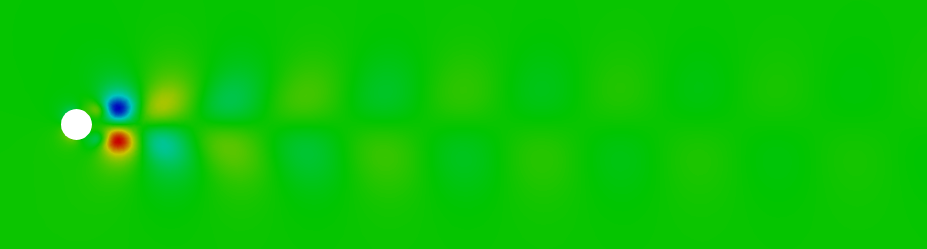}\hspace*{1cm}%
  \includegraphics[width=\linewidth]{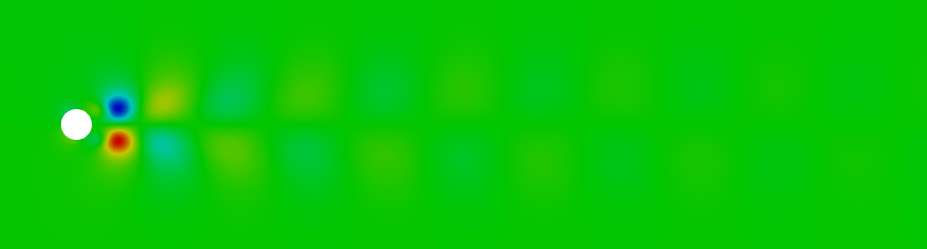}
  \caption{mode 2}
  \label{fig:sub1}
\end{subfigure}
\par\bigskip
\begin{subfigure}{.5\textwidth}
  \centering
  \hspace*{-4cm}\includegraphics[width=\linewidth]{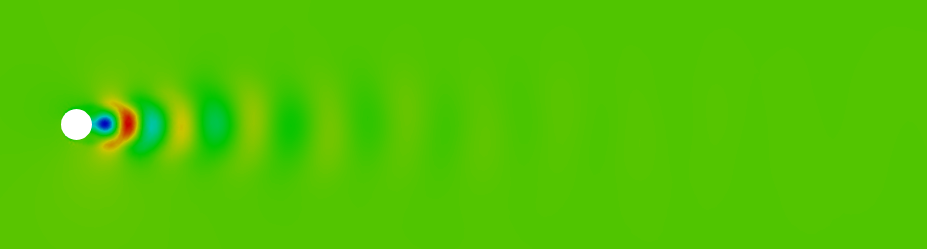}\hspace*{1cm}%
  \includegraphics[width=\linewidth]{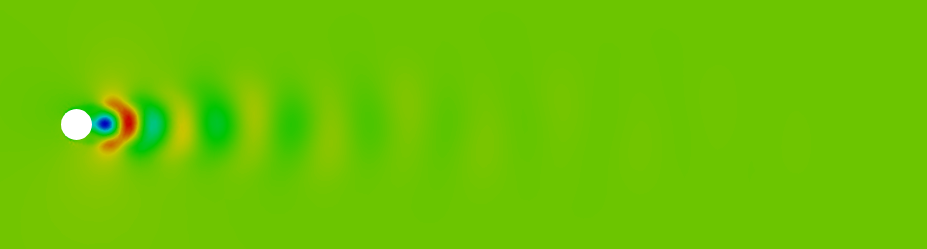}
  \caption{mode 3}
  \label{fig:sub1}
\end{subfigure}
\par\bigskip
\begin{subfigure}{.5\textwidth}
  \centering
  \hspace*{-4cm}\includegraphics[width=\linewidth]{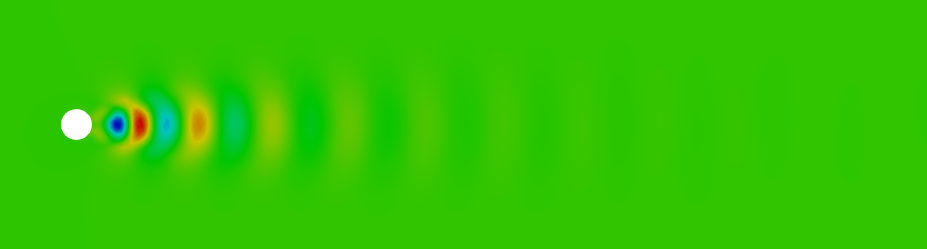}\hspace*{1cm}%
  \includegraphics[width=\linewidth]{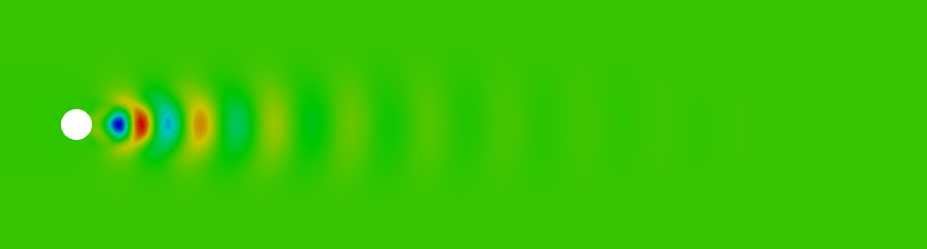}
  \caption{mode 4}
  \label{fig:sub1}
\end{subfigure}
\caption{POD (left) and calibrated ITSGM (right) pressure modes associated to the untrained Reynolds number $Re = 195$.}
\label{fig:pressure_space_modes_ITSGM_195}
\end{figure}

\begin{figure}[hbtp!]
%\centering 
\hspace*{-2cm}
\includegraphics[width=1.2\linewidth]{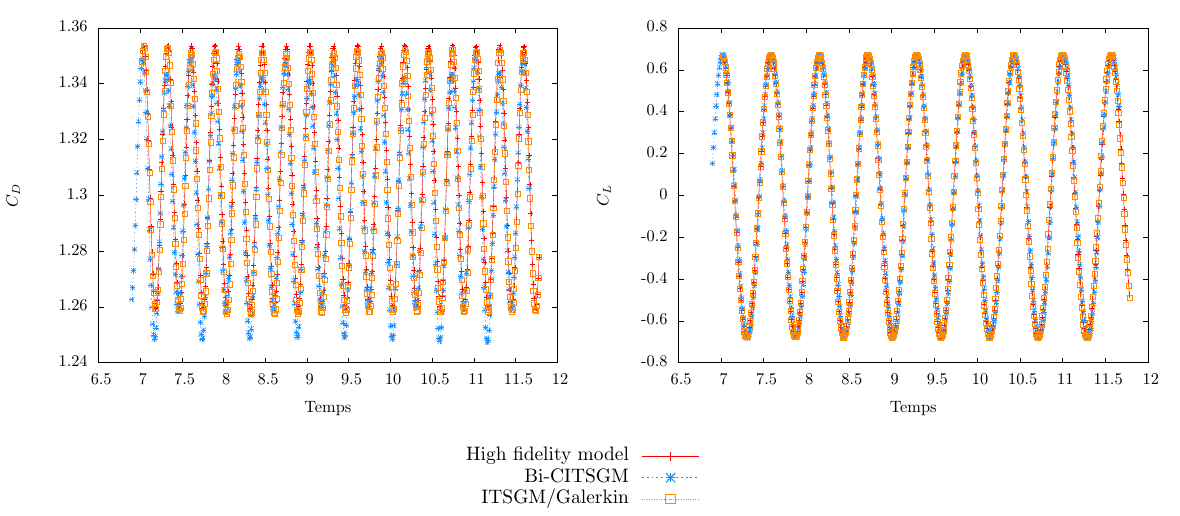}
\caption{Lift and drag coefficients $C_D$ and $C_L$ obtained by Bi-CITSGM and ITSGM/Galerkin with those obtained by the high fidelity model at the untrained Reynolds number $Re = 195$.}
\label{fig:results_CD_CL_195}
\end{figure}

\begin{figure}[hbtp!]
%\centering 
\hspace*{-2cm}
\includegraphics[width=1.2\linewidth]{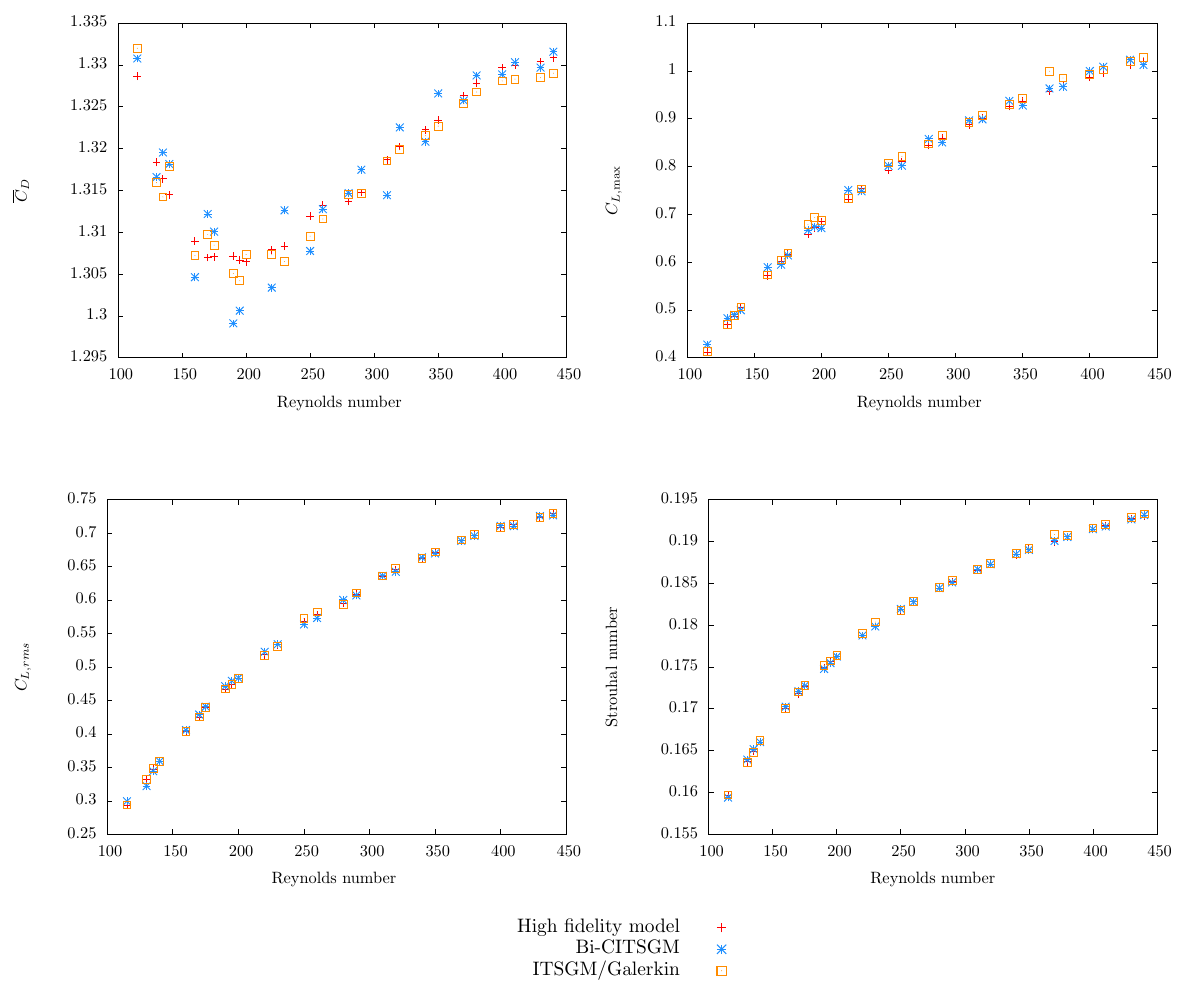}
\caption{Aerodynamics coefficient and Strouhal number at different untrained Reynolds number values.}
\label{fig:hydrodynamics_coeff_Re_range}
\end{figure}

\begin{figure}[hbtp!]
%\centering 
\hspace*{-2cm}
\includegraphics[width=1.2\linewidth]{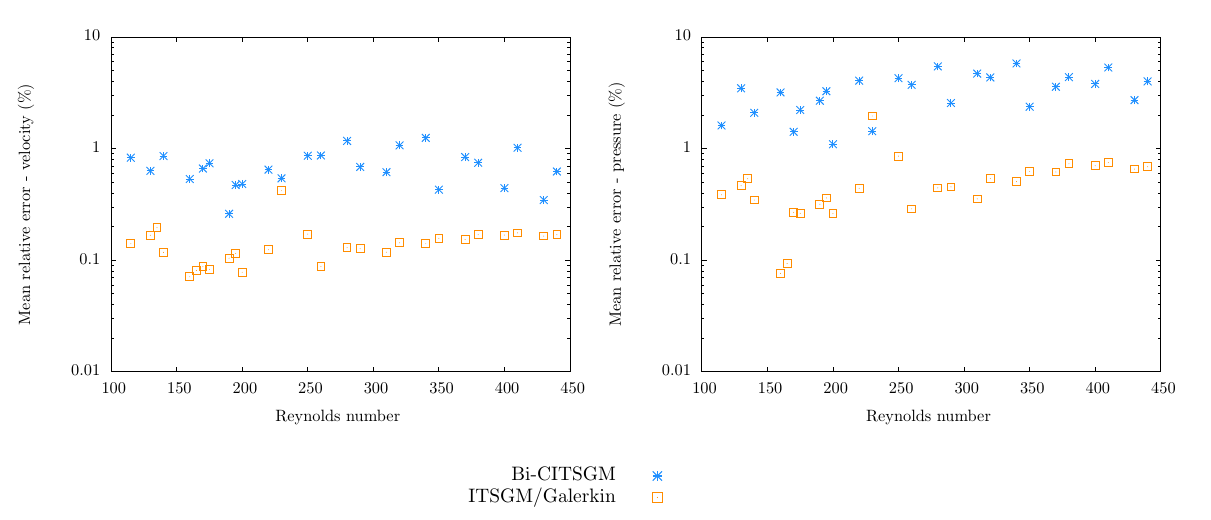}
\caption{Mean relative Errors at different untrained Reynolds number values.}
\label{fig:Errors_Re_range}
\end{figure}

%For different values of the untrained Reynolds number $\widetilde{Re}$, using the recovered SVD form of the trained snapshots $\SnapMatrix_{u,i}$ and $\SnapMatrix_{p,i}$ $i=1,\dots,\nbrBases$, the NIMR method can be applied to find an approximation of $\widetilde{\SnapMatrix}_u$ and $\widetilde{\SnapMatrix}_p$. In this example, $3$ neighbour bases are selected in the application of the ITSGM method and the weights inverse powers $m$ and $l$ are such as $m=l=5/2$. Moreover, Lagrange method is considered for interpolation of initial geodesic velocities in the tangent subspace of the Grassmann manifold and cubic spline for singular values interpolation. 
For different untrained Reynolds number values, the $\nonintrusivemethod$ method is applied to find an approximation of the solutions snapshots of both velocity and pressure. Cubic spline method was used to interpolate the singular values. Whilst ITSGM is used for bases interpolation. In this example, Lagrange interpolation is used in the tangent space of the Grassmann manifold, where only initial geodesic velocities associated to the $3$ closest trained Reynolds number values are interpolated. The inverse weights powers $m$ and $l$ are chosen such as $m=l=3$.
In a first attempt to validate the Bi-CITSGM approach, the untrained  Reynolds number $Re = 195$, was considered. The first $4$ POD and Bi-CITSGM basis functions are represented in figures \ref{fig:velocity_space_modes_ITSGM_195} and \ref{fig:pressure_space_modes_ITSGM_195}. It can be seen that the Bi-CITSGM spatial modes are almost identical 
%succeed to fit very well with 
to the original POD modes.
%Table \ref{Tab:NIMR_Re195} lists the aerodynamics coefficients and errors with respect to the high fidelity solution of POD, ITSGM/ROM and NIMR solutions. It can be seen that 
%
%
The quality of the Bi-CITSGM temporal modes can be checked out from inspecting the hydrodynamics coefficients. Figure \ref{fig:results_CD_CL_195} shows a good match between lift and drag coefficients $C_L$ and $C_D$, obtained by Bi-CITSGM and ITSGM/Galerkin with the coefficients obtained by the high fidelity model.
%Visually, a slight difference can be noticed. In order to see this difference, table \ref{Tab:NIMR_Re195} lists the corresponding  aerodynamics coefficients.
This agreement of results is shown in the mean hydrodynamics coefficients listed in table \ref{Tab:NIMR_Re195}. 
%As can be seen from the table, the Bi-CITSGM allowed a good prediction of these coefficients. 
Let's denote by $\overline{\varepsilon}^{\%}$ the mean relative error given by
$$\overline{\varepsilon}^{\%} = 100 \times \left. \left( \int_0^T||f- \widetilde{f} ||_{L^2(\Omega)}^2  \, dt \right)^{\frac{1}{2}} \middle/ \left( \int_0^T|| f ||_{L^2(\Omega)}^2  \, dt \right)^{\frac{1}{2}} \right.$$
where $\widetilde{f}$ is an approximation of $f$.
It is confirmed by inspecting this error for the untrained Reynolds number $Re = 195$, that the Bi-CITSGM presented a low error of about $0.4\%$ for velocity and $3.2\%$ for pressure. 
%
%
%In addition, figure \ref{fig:results_CD_CL_195} shows a good fitting between lift and drag coefficients $C_L$ and $C_D$, obtained by NIMR and the high fidelity model.
%  
%
%
\\
Now, in order to demonstrate the robustness of this method, the same previous study was performed for several untrained Reynolds number values in the interval $]90, 450[$. In figure \ref{fig:hydrodynamics_coeff_Re_range}, the aerodynamics coefficients $\overline{C}_D$,  $C_{L,max}$, $C_{L,rms}$ and Strouhal number $S_t$ obtained by ITSGM/Galerkin and Bi-CITSGM are compared to those obtained by the high fidelity model \footnote{The shape of the coefficient $\overline{C}_D$ in figure \ref{fig:hydrodynamics_coeff_Re_range} shows a convex behaviour with a minimum included between $Re=160$ and $Re=230$. The same behaviour was observed by Fabiane et al. \cite{Fiabane-2011}.}. For $C_{L,max}$, $C_{L,rms}$ and $S_t$, a good agreement between the high fidelity results and Bi-CITSGM method is observed, whilst slight differences are observed on the mean drag coefficient $\overline{C}_D$. In terms of mean relative errors, it can be seen from figure \ref{fig:Errors_Re_range} that for all the tested untrained Reynolds number values, the mean error was less than $1.5\%$ for velocity and $6\%$ for pressure.
Finally, table \ref{Tab:Temps_CPU_interp} reports that the proposed method Bi-CITSGM  performed in real time (less than one second). Concretely, the non-intrusive approach Bi-CITSGM is about $10$ times faster than the ITSGM/Galerkin approach.
%
%Finally, table \ref{Tab:Temps_CPU_interp} reports that the Bi-CITSGM is computationally efficient compared to ITSGM/Galerkin with a ratio of about $10$ times. 
The reason is that ITSGM/Galerkin requires more computational operations involved essentially in the Galerkin projection process, while only few computational operations are needed for the Bi-CITSGM method.
%
%
%
%
%
%
%
%
%
%
%Le pourcentage d'erreur moyen $\overline{\varepsilon}_{\%}$ entre la solution du modele complet $\vect{\stateU}$ et son approximation $\vect{\stateU}_{\nbrModes_u}$ est mesure en utilisant la relation
%\begin{equation*}
%\overline{\varepsilon}_{\vect{\stateU}}^{\%} = 100\times \myfrac{ \left( \myint{0}{T}\normLtwo{\vect{\stateU}-\vect{\stateU}_{\nbrModes}}^2\, dt \right)^{\frac{1}{2}}}{ \left(\myint{0}{T}\normLtwo{\vect{\stateU}}^2\, dt\right)^{\frac{1}{2}}} 
%\end{equation*}
%Les mêmes relations sont utilises pour mesurer l'erreur d'approximation de la pression $p$. 
%
%
%
%
%
%
%
%
%
%% BEGIN TABLE ------------------------------------------------------------------------------------------------------------------------------
\begin{table}[hbtp!]
\small
 \centering
% \begin{tabular}{ |p{2.5cm}|p{1.5cm}| p{2.7cm}| p{2.7cm}|p{1.2cm}| }
\begin{tabular}{ |c|c|c|c|c|c|c| }
 \hline
Methode &$\overline{C}_D$&  $C_{L,max}$    & $C_{L,rms}$ & $S_t$ & $\overline{\varepsilon}_{\vect{\stateU}}^{\%}$ & $\overline{\varepsilon}_p^{\%}$ \\
\hline
%----------------------------------------------------------------------------------------------------
High fidelity model	&	$1.30$	&	$0.67$	&	$0.47$	&	$0.17$	&	-	&	- \\ 
%POD	&	$1.307$	&	$0.672$	&	$0.474$	&	$0.176$	&	$0.03426\%$	&	$0.03631\%$\\ 
%ROM	&	$1.306$	&	$0.694$	&	$0.474$	&	$0.176$	&	$0.09689\%$	&	$0.31798\%$\\ 
ITSGM/Galerkin	&	$1.30$	&	$0.69$	&	$0.47$	&	$0.17$	&	$0.1\%$	&	$0.3\%$\\ 
Bi-CITSGM	&	$1.30$	&	$0.67$	&	$0.48$	&	$0.17$	&	$0.4\%$	&	$3.2\%$\\ 
%----------------------------------------------------------------------------------------------------
\hline
\end{tabular}
\caption{Aerodynamics coefficients and mean relative at the untrained Reynolds number $Re = 195$. $\overline{\varepsilon}_{\vect{\stateU}}^{\%}$ and $\overline{\varepsilon}_p^{\%}$ are respectively velocity and pressure mean relative errors.}
\label{Tab:NIMR_Re195}
\end{table}
% END TABLE  ------------------------------------------------------------------------------------------------------------------------------
%
%
%
%
%
% END FIGURE -------------------------------------------------------------------------------------------------------------------------------
%
%
%
%
%
%
%\begin{figure}[hbtp!]
%\hspace*{-0.8 in}
%\begin{subfigure}{.6\textwidth}
%  \centering
%  \resizebox{\textwidth}{!}{\input{./FIGURES/PART3/NIMR/mean_aerodynamics_coeffs/CD_mean}}
%  %\caption{Coefficient de traînee moyen $\overline{C}_D$}
%  \label{fig:sub1}
%\end{subfigure}%
%\begin{subfigure}{.6\textwidth}
%  \centering
%  \resizebox{\textwidth}{!}{\input{./FIGURES/PART3/NIMR/mean_aerodynamics_coeffs/CL_max}}
%  %\caption{Coefficient de portance maximum $C_{L,max}$}
%  \label{fig:sub2}
%\end{subfigure}
%\par\bigskip\centering
%\hspace*{-0.8 in} % force a bit of vertical whitespace
%\begin{subfigure}{.6\textwidth}
%  \centering
%  \resizebox{\textwidth}{!}{\input{./FIGURES/PART3/NIMR/mean_aerodynamics_coeffs/Crms}}
%  %\caption{Crms}
%  \label{fig:sub3}
%\end{subfigure}%
%\begin{subfigure}{.6\textwidth}
%  \centering
%  \resizebox{\textwidth}{!}{\input{./FIGURES/PART3/NIMR/mean_aerodynamics_coeffs/strouhal}}
%  %\caption{Nombre de Strouhal}
%  \label{fig:sub3}
%\end{subfigure}
%\\ 
%\centering
%\hspace*{-1cm}\adjustbox{width= 0.6\linewidth,trim= 1cm 4cm 1cm 4cm}{\input{./FIGURES/PART3/NIMR/mean_aerodynamics_coeffs/LEGEND}}
%\caption{Aerodynamics coefficient and Strouhal number at different untrained Reynolds number values.}
%\label{fig:hydrodynamics_coeff_Re_range}
%\end{figure}

% BEGIN FIGURE -----------------------------------------------------------------------------------------------------------------------------
% END FIGURE -------------------------------------------------------------------------------------------------------------------------------
%
%
%
%
%%% BEGIN TABLE --------------------------------------------------------------------------------------------------------------------------
\begin{table}[hbtp!]
\small
 \centering
\begin{tabular}{ |c|c| }
 \hline
ITSGM/Galerkin & Bi-CITSGM    \\
\hline
%----------------------------------------------------------------------------------------------------
$9.5$ sec	&	$0.8$ sec	\\
%----------------------------------------------------------------------------------------------------
\hline
\end{tabular}
\caption{CPU computational time}
\label{Tab:Temps_CPU_interp}
\end{table}
% END TABLE  ----------------------------------------------------------------------------------------------------------------------------

 %%%%%%%%%%%%%%%%%%%%%%%%%%%%%%%%%%%%%%%%%%%%%%%%%%%%%%%%%%%%%%%%%%%%%%%%%%%%%%%%%%%%%%%%%%%%%%%%%%%%%%%%%%%%%%%%%%%%%%%%%%%%%%%%%%%%%%%%%%%%%%%%%%%%
\section{Conclusions}
In this article, we proposed a non intrusive method Bi-CITSGM for model reduction of parametrized non linear time-dependent physical problems. The steps of the method are to represent a sampling of parametrized solutions snapshots as an ensemble of POD decompositions, approximate singular values (square roots of POD eigenvalues) using a spline cubic interpolation method, interpolate the spatial and temporal bases (left and right singular vectors) using ITSGM, and finally calibrate the result by two orthogonal matrices obtained as  analytical solutions of two constrained optimization problems. The ability of this approach to quickly and accurately reproduce the solution snapshots for new untrained parameters was numerically investigated in the problem case of the flow past a circular cylinder, where the parameter of interpolation was considered as the Reynolds number (ranged between $90$ and $450$). The accuracy of the velocity and pressure solution snapshots obtained for different untrained Reynolds numbers reproduced by the Bi-CITSGM was checked.  The mean relative error was less than $1.5\%$ for velocity and $6\%$ for the pressure.
The computational time needed for the Bi-CITSGM (less than $1$ second) is reduced by one order of magnitude with respect to the time needed for ITSGM/Galerkin, and by several orders of magnitude with respect to the time needed for the high fidelity model. Future work includes applying the Bi-CITSGM approach on smolyak grids to more realistic scenarios (experimental data), and eventually use it with a genetic algorithm to solve inverse problems.

 \section*{Acknowledgement}
 This material is based upon work financially supported by the Nouvelle-Aquitaine region and CPER/FEDER b\^atiment durable.

\appendix
\section{Review of POD}
Let $\SnapMatrix\in \mbb{R}^{\dimVecorsSnap\times \nbrSnap}$ be the snapshot matrix represented in the low dimensional subspace $\txt{span}(\Xbasis{})$ obtained by the POD method. Here, the POD basis can be chosen to be optimal with respect to any inner product. Assume that optimality is realized with respect to an inner product \footnote{Different choices can be made for the inner product of optimality of the POD basis. For example, if the POD basis is to be optimal with respect to the $L^2$ inner product, the operator $\mcal{A}$ is the mass matrix. The simplest choice is the Euclidean inner product for which $\mcal{A}$ is the identity matrix.} for which the associated linear operator is denoted $\mcal{A}$. The POD method consists of the following steps
\begin{itemize}
\item[\hspace{10pt}]
\begin{itemize}
\item[\textit{step 1}] build the correlation matrix $C$ as $C = \SnapMatrix^T \mcal{A} \SnapMatrix$
\item[\textit{step 2}] solve the eigenvalue problem $C \Tbasis{} = \Tbasis{} \lambda $
\item[\textit{step 3}] calculate the POD left singular vectors as $\Xbasis{} = \SnapMatrix \Tbasis{} \lambda^{-\frac{1}{2}}$
%\item[\textit{step 4}] calculate the POD right singular vectors as $\Tbasis{} = $
\end{itemize}
\end{itemize}
Let $\Sigma = \lambda^{\frac{1}{2}}$, the snapshots matrix $\SnapMatrix$ can be written as
\begin{equation*}
\SnapMatrix = \Xbasis{} \Sigma \Tbasis{}^T
\end{equation*}
In applications, the POD basis is truncated to an order $\nbrModes<\nbrSnap$, where only the first modes corresponding the the first $\nbrModes$ significant eigenvalues of $C$ are considered. This truncated basis is sufficient to represent the most of the information contained in $\SnapMatrix$.
%where $\Sbase{i}\in \mbb{R}^{\dimVecorsSnap\times \nbrModes}$ and $\Tbase{i} \in \mbb{R}^{\nbrSnap \times \nbrModes}$ verify
%\begin{equation}
%\integOmega{\Smode{i}{j} \Smode{i}{k}} = \delta_{jk} \hspace*{0.5cm} \txt{and} \hspace*{0.5cm} \Tmode{i}{j}^T \Tmode{i}{k} = \tau_{jk} \delta_{ij}
%\end{equation}
%Let $\Mat{U}_i$ the basis obtained by applying the Gram Schmidt Orthonormalization to $\Sbase{i}$, it is possible to recover the SVD form of $\SnapMatrix_i$ as follows
%\begin{equation*}
%\SnapMatrix_i \approx \Mat{U}_{i} \Mat{\Sigma}_i \Mat{V}_{i}^T
%\end{equation*}
%where $\Mat{U}_{i} = \Sbase{i}\Tbase{i}^T\Mat{V}_{i} \Mat{\Sigma}_i^{-\frac{1}{2}}$ such that $\Mat{U}_i^T \Mat{U}_i = I_{\nbrModes}$ and $\Mat{\Sigma}_i$ and $\Mat{V}_{i}$ verify the following eigendecomposition
%\begin{equation}
%\Tbase{i}\Sbase{i}^T\Sbase{i} \Tbase{i}^T \Mat{V}_{i} =  \Mat{\Sigma}_i \Mat{V}_{i}
%\end{equation}
\section{POD reduced order model of the flow past a cylinder}\label{Appendix.Cylinder}
By plugging expressions \eqref{EQ : subspace approx POD u and p} into equations \eqref{EQ: Navier_Stokes_cylinder} and imposing the orthogonality condition of the residual on the basis functions $\Smode{u}{i}$ and on the gradient of $\Smode{p}{l}$, the reduced order model associated to the problem of flow past a cylinder writes
\begin{equation}\label{Eq : ROM CYLINDER}
\begin{cases}
\somme{j}{1}{\nbrModesU}\Mat{M}^{(u)}_{ij} \derivee{\Tmode{u}{j}}{t} + \somme{j}{1}{\nbrModesU}\left[ \myfrac{1}{Re} \Mat{R}^{(u)}_{ij} +  \Mat{\overline{C}}^{(u)}_{ij} \right] \Tmode{u}{j} +  \somme{j}{1}{\nbrModesU}\somme{k}{1}{\nbrModesU}\Mat{C}^{(u)}_{ijk}  \Tmode{u}{j}\Tmode{u}{k} + \somme{l}{1}{\nbrModesP} \Mat{K}^{(u)}_{il} \Tmode{p}{l}= \widetilde{\bm{F}}^{(u)}_{i}, \ 
\\
\somme{j}{1}{\nbrModesU} \Mat{M}^{(p)}_{mj} \derivee{\Tmode{u}{j}}{t} + \somme{j}{1}{\nbrModesU} \left[ \myfrac{1}{Re}\Mat{R}^{(p)}_{mj}+\Mat{\overline{C}}^{(p)}_{mj}\right]   \Tmode{u}{j} + \somme{j}{1}{\nbrModesU}\somme{k}{1}{\nbrModesU} \Mat{C}^{(p)}_{mjk} \Tmode{u}{j}\Tmode{u}{k}+ \somme{l}{1}{\nbrModesP} \Mat{K}^{(p)}_{ml} \Tmode{p}{l}= \widetilde{\bm{F}}^{(p)}_{m}
\\
\somme{j}{1}{\nbrModesU}\Mat{M}^{(u)}_{ij} \Tmode{u}{j}(0) = \int_{\Omega} \vect{u_0}  \Smode{u}{i}\,dx
\\ 
\forall i = 1,\cdots, \nbrModesU, \ \ \ \forall m = 1,\cdots, \nbrModesP
\end{cases}
\end{equation}
where 
\begin{align*}
\Mat{M}^{(u)}_{ij} &= \integOmega{\Smode{u}{j} \Smode{u}{i}}
\hspace*{1cm}
\Mat{R}^{(u)}_{ij}  = \integOmega{\nabla \Smode{u}{j} : \nabla\Smode{u}{i}} 
\hspace*{1cm} 
\Mat{K}^{(u)}_{il}  = \int_{\Gamma} \Smode{p}{l} \Smode{u}{i}\cdot \vect{n} \,d\sigma
\\
\Mat{\overline{C}}^{(u)}_{ij} &= \integOmega{(\overline{\statevecU}\cdot\nabla) \Smode{u}{j} \cdot \Smode{u}{i}} + \integOmega{(\Smode{u}{j}\cdot\nabla) \overline{\statevecU} \cdot \Smode{u}{i}}\\
\Mat{C}^{(u)}_{ijk} &= \integOmega{(\Smode{u}{j}\cdot\nabla) \Smode{u}{k} \cdot \Smode{u}{i}} 
\hspace*{2cm}
\widetilde{\bm{F}}^{(u)}_{i} = \integOmega{ \left(\myfrac{1}{Re} \Delta \overline{\statevecU} - \overline{\statevecU} \cdot\nabla \overline{\statevecU} - \nabla \overline{p}\right) \, \Smode{u}{i}}
\\
\Mat{M}^{(p)}_{mj} &= \integOmega{\Smode{u}{j} \nabla\Smode{p}{m}}
\hspace*{1cm}
\Mat{R}^{(p)}_{mj}  = -\integOmega{\Delta \Smode{u}{j} : \nabla\Smode{p}{m}} 
\hspace*{1cm} 
\Mat{K}^{(p)}_{ml}  = \int_{\Omega} \nabla\Smode{p}{l} \nabla\Smode{p}{m} \,dx
\\
\Mat{\overline{C}}^{(u)}_{mj} &= \integOmega{(\overline{\statevecU}\cdot\nabla) \Smode{u}{j} \cdot \nabla\Smode{p}{m}} + \integOmega{(\Smode{u}{j}\cdot\nabla) \overline{\statevecU} \cdot \nabla\Smode{p}{m}}\\
\Mat{C}^{(p)}_{mjk} &= \integOmega{(\Smode{u}{j}\cdot\nabla) \Smode{u}{k} \cdot \nabla\Smode{p}{m}} 
\hspace*{2cm}
\widetilde{\bm{F}}^{(p)}_{m} = \integOmega{ \left(\myfrac{1}{Re} \Delta \overline{\statevecU} - \overline{\statevecU} \cdot\nabla \overline{\statevecU} - \nabla \overline{p}\right) \, \nabla\Smode{p}{m}}
\end{align*}
More details about this reduced order model can be found in \cite{Tallet-Allery-com-2015}.
 %%%%%%%%%%%%%%%%%%%%%%%%%%%%%%%%%%%%%%%%%%%%%%%%%%%%%%%%%%%%%%%%%%%%%%%%%%%%%%%%%%%%%%%%%%%%%%%%%%%%%%%%%%%%%%%%%%%%%%%%%%%%%%%%%%%%%%%%%%%%%%%%%%%%

%\section*{\large References}
\bibliographystyle{ieeetr}
\bibliography{./BIBLIO}
\end{document}